\documentclass{article}

\usepackage{graphicx}
\usepackage[german,english]{babel}
\usepackage{amsmath}
\usepackage{amssymb}
\usepackage{amsfonts}
\usepackage{bbm}
\usepackage{url}
\usepackage{booktabs}
\usepackage{tabularx}
\usepackage{epsfig,subfigure}
\usepackage{longtable}
\usepackage{lscape}
\usepackage{psfrag}

\newtheorem{defn}{Definition}

\newtheorem{thm}[defn]{Theorem}

\newtheorem{conj}[defn]{Conjecture}

\makeatletter

\title{Triangulated Manifolds with Few Vertices: Geometric $3$-Manifolds}

\author{\Large Frank H.~Lutz}
\date{}

\begin{document}

\selectlanguage{english}

\maketitle

\bigskip
\bigskip
\bigskip
\bigskip
\bigskip
\bigskip


The understanding and classification of (compact) $3$-dimensional manifolds
(without boundary) is with no doubt one of the most prominent tasks 
in topology ever since Poincar\'e's fundamental work \cite{Poincare1904} 
on \guillemotleft l'analysis situs\guillemotright\ \linebreak
appeared in~1904.

There are various ways for constructing $3$-manifolds, some 
of which that are general enough to yield all $3$-manifolds 
(orientable or nonorientable) and some that produce only 
particular types or classes of examples.

According to Moise \cite{Moise1952}, all $3$-manifolds can be
triangulated. This implies that there are only countably many 
distinct combinatorial (and therefore at most so many different topological) types 
that result from gluing together tetrahedra. Another 
way to obtain $3$-manifolds is by starting with a solid $3$-dimensional polyhedron 
for which surface faces are pairwise identified
(see, e.g., Seifert~\cite{Seifert1931} and Weber and Seifert \cite{WeberSeifert1933}).
Both approaches are rather general and, on the first sight, 
do not give much control on the kind of manifold we can expect
as an outcome. However, if we want to determine the topological type
of some given triangulated $3$-manifold, then small or minimal
triangulations play an important role as reference objects 
for recognition heuristics based on bistellar flips; see \cite{BjoernerLutz2000}.

Further frequently used methods to get all orientable $3$-manifolds
are via composition of two handle bodies of the same genus
(\emph{Heegaard splitting}; \linebreak
see Heegaard \cite{Heegard1916}, Seifert and Threlfall \cite[\S 63]{SeifertThrelfall1980}, 
Hempel \cite[Ch.~2]{Hempel1976} and Stillwell \cite[Ch.~8.3]{Stillwell1980}),
via surgery on a link in the $3$-sphere $S^3$ 
(cf.\  Dehn \cite{Dehn1910}, Seifert and Threlfall \cite[\S 65]{SeifertThrelfall1934},
Lickorish \cite{Lickorish1962}, and Wallace \cite{Wallace1960})
or via three-fold branched coverings of $S^3$ over a knot 
(see Alexander \cite{Alexander1919-1920},
Hilden \cite{Hilden1974}, Montesinos \cite{Montesinos1976},
and Izmestiev and Joswig \cite{IzmestievJoswig2001pre}).

In this article, we explicitly construct small triangulations 
for a number of well-known $3$-dimensional manifolds
and give a brief outline of some aspects of the underlying 
theory of $3$-manifolds and its historical development.


\section{Geometrization of $3$-Manifolds}

The intensive study of non-Euclidean geometries in the late 19th century,
in particular by Klein \cite{Klein1890}, led, as one consequence, 
to the famous \emph{Clifford-Klein space form problem}.
This problem, posed by Killing \cite{Killing1891} in 1891,
asks for a classification of all (Riemannian) manifolds of constant curvature 
that arise as quotients $X^d/G$, where $X^d$ is either the $d$-dimensional sphere $S^d$, 
the  $d$-dimensional Euclidean space $E^d$, or the $d$-dimensional
hyperbolic space~$H^d$, and $G$ is a discrete group of isometries which acts freely on $X^d$.
The resulting manifolds are called \emph{spherical}, \emph{flat}, or \emph{hyperbolic},
respectively. (This modern formulation of the Clifford-Klein space form problem
is by Hopf \cite{Hopf1926}.)

The three-dimensional spherical space forms were classified by
Hopf \cite{Hopf1926} and by
Threlfall and Seifert (\cite{ThrelfallSeifert1931}, \cite{ThrelfallSeifert1933}).
The classification of compact flat $3$-manifolds 
is due to Nowacki \cite{Nowacki1934} up to homeomorphisms,
due to Hantzsche and Wendt \cite{HantzscheWendt1935} in the affine case, 
and was generalized by Wolf \cite[Sec.\ 3.5]{Wolf1967} to the isometric and
noncompact case.
Hyperbolic $3$-manifolds turned out to be much harder to understand.
At the beginning, it was not even clear whether there are 
closed hyperbolic $3$-manifolds at all.
A first series of examples was constructed by L\"obell in 1931 \cite{Loebell1931},
and in 1933 Weber and Seifert \cite{WeberSeifert1933} 
presented their hyperbolic dodecahedral space.
The picture completely changed when Thurston \cite{Thurston1982}
discovered that almost every prime $3$-manifold is hyperbolic.
For surveys on hyperbolic $3$-manifolds see Thurston (\cite{Thurston1997}, \cite{Thurston2002}),
Benedetti and Petronio \cite{BenedettiPetronio1992}, 
and McMullen \cite{CMcMullen1992}. 
Milnor \cite{Milnor1982} is an exposition 
of the development of hyperbolic geometry. Commented English translations of 
classical papers by Beltrami, Klein, and Poincar\'e on hyperbolic
geometry can be found in Stillwell \cite{Stillwell1996}. 

In dimension two, every surface $M^2$ is a Clifford-Klein
manifold. If the Euler characteristic $\chi (M^2)$ is positive, 
then $M^2$ is spherical, which gives the $2$-sphere $S^2$ 
and the real projective plane ${\mathbb R}{\bf P}^{\,2}$.
The only two flat $2$-manifolds with $\chi =0$ are the $2$-torus $T^2$ and
the Klein bottle $K$. All other surfaces are hyperbolic and have negative 
Euler characteristic.

In dimension three, 
there are, besides the Clifford-Klein spaces, 
manifolds that are \emph{homogeneous} but not \emph{isotropic} 
(i.e., they look the same at every point but differently 
when we look in different directions).
The sphere product $S^2\!\times\!S^1$ is an example for such a space: 
This manifold can be modeled as a quotient 
of the simply connected space $S^2\times{\mathbb R}^1$.

In a unifying approach for a geometric description
of all $3$-manifolds, \linebreak Thurston \cite{Thurston1982} introduced 
the concept of \emph{model geometries} as
potential building blocks. For surveys on the geometrization
of $3$-manifolds see the electronic edition \cite{Thurston2002} 
of the 1980 lecture notes and the book \cite{Thurston1997} of Thurston
as well as the comprehensive article of Scott \cite{Scott1983}.

\begin{defn} {\rm (Thurston \cite[3.8.1]{Thurston1997})}
A \textbf{model geometry} $(X,G)$ consists
of a simply connected (and connected) manifold $X$
together with a transitive Lie group $G$ of diffeomorphisms of $X$,
such that $G$ has compact point stabilizers and is maximal with this
property. Moreover it is required that there is at least
one compact manifold $M$ modeled on $(X,G)$ as a quotient
$X/H$ with respect to a discrete fixed point free subgroup $H\subset G$.
\end{defn}
Any model geometry carries a homogeneous $G$-invariant metric.
For dimension two the model spaces $X$ are precisely the
spaces of constant curvature $S^2$, $E^2$, and $H^2$.
\begin{thm} {\rm (Thurston \cite[3.8.4]{Thurston1997})}
There are eight $3$-dimensional model geometries $(X,G)$,
namely the isotropic ones $S^3$, $E^3$, and $H^3$,
and the anisotropic ones\, $S^2\times{\mathbb R}^1$, $H^2\times{\mathbb R}^1$,
$\widetilde{SL}(2,{\mathbb R})$, ${\rm Nil}$, and\, ${\rm Sol}$.
\end{thm}
If a closed $3$-manifold can be modeled on one of the eight geometries,
then the respective geometry is unique. 
\begin{conj} {\rm (Thurston's Geometrization Conjecture \cite{Thurston1982})}
Every compact $3$-manifold can be decomposed canonically into geometric pieces.
\end{conj}
The decomposition is as follows. In a first step, the manifold is split
into a connected sum of (finitely many) prime factors; 
see Kneser \cite{Kneser1929}, Milnor \cite{Milnor1962}, Haken \cite{Haken1961a}, 
and Jaco and Tollefson \cite{JacoTollefson1995}.
In a second step, each prime factor is then cut along a finite family
of certain nontrivially embedded tori. This second step was formulated
by Johannson \cite{Johannson1979} and Jaco and Shalen \cite{JacoShalen1979}
(see also Thurston \cite{Thurston1982}, Scott \cite{Scott1983}
and Neumann and Swarup \cite{NeumannSwarup1997})
and leaves as the open part of the conjecture to show that each 
resulting piece can be modeled on one of the eight geometries.

Thurston's Geometrization Conjecture holds in many cases.
For example, all compact, prime $3$-manifolds with nonempty boundary 
and all Haken mani\-folds are geometric \cite{Thurston1986}.
Recent progress towards the Geometrization Conjecture has been achieved
by Perelman (\cite{Perelman2002pre}, \cite{Perelman2003bpre}, \cite{Perelman2003apre}),
based on work by Hamilton (\cite{Hamilton1982}, \cite{Hamilton1999}) 
on the Ricci flow on three-manifolds.

All eight model geometries are Lie groups, and, with the exception
of $H^3$, all are unimodular. Compact quotients 
of the spaces $S^3$, $E^3$, $\widetilde{SL}(2,{\mathbb R})$, 
${\rm Nil}$, and ${\rm Sol}$ that arise from discrete 
subgroups acting by left (or right) multiplication
were classified by Raymond and Vasquez \cite{RaymondVasquez1981},
based on partial results by Auslander, Green, 
and Hahn \cite[Ch.\ III]{AuslanderGreenHahn1963}.

Six of the geometries, $S^2\times{\mathbb R}^1$, $S^3$, $E^3$, ${\rm Nil}$,
$H^2\times{\mathbb R}^1$, and $\widetilde{SL}(2,{\mathbb R})$,
are Seifert manifolds and will be discussed in the next sections.


\section{Seifert Manifolds}

In 1933, Seifert (see \cite{Seifert1933} and the 
English translation \cite{SeifertThrelfall1980}) described
a large class of $3$-manifolds, which revolutionized 
$3$-manifold theory at that time. 
In the following we give a summary 
of the main definitions and theorems following Seifert \cite{Seifert1933}
and Orlik \cite{Orlik1972}.

\begin{defn}
A $3$-manifold is a \textbf{Seifert manifold}
if it can be decomposed into a disjoint union
of simple closed curves so that every such fiber 
has a (closed) fibered solid torus as a tubular neighborhood.
\end{defn}

Any \emph{fibered solid torus} is obtained by starting with
an upright cylinder for which the top and the bottom side
are identified by a twist of a rational angle $2\pi \frac{\nu}{\mu}$, 
where $\mu$ and $\nu$ are coprime integers. 
If the cylinder is fibered into straight line segments parallel to its center
axis, then the resulting torus is fibered into closed circles.
Replacing $\nu$ by $\nu+k\mu$ or by $-\nu$ yields, 
up to a fiber preserving homeomorphism, the same 
fibered solid torus. Therefore, without loss of generality, 
it can be assumed that $\mu>0$ and $0\leq\nu\leq\frac{1}{2}\mu$.

\begin{figure}
\begin{center}
\footnotesize
\psfrag{1}{1}
\psfrag{2}{2}
\psfrag{3}{3}
\psfrag{4}{4}
\psfrag{5}{5}
\includegraphics[width=.46\linewidth]{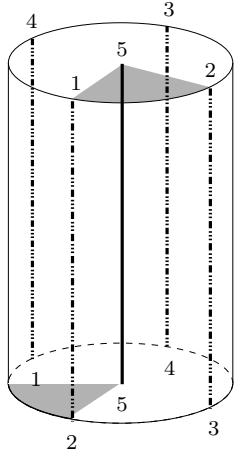}
\end{center}
\caption{Fibered solid torus with $\mu=4$ and $\nu=1$.}
\label{fig:cylinder1}
\end{figure}

Figure~\ref{fig:cylinder1} displays a fibered solid torus with $\mu=4$
and $\nu=1$, i.e., the top and the bottom of the cylinder are identified
by a twist of $\pi/2$. The dashed fiber is composed of the
four line segments 1--2, 2--3, 3--4, 4--1 on the toroidal boundary surface. 
Every other fiber of the solid torus is composed of four line segments
as well, with the exception of the center fiber, which consists of one
line segment only.

Two fibered solid tori can be mapped onto each other by a 
fiber preserving homeomorphism if and only if
the respective parameters $\mu$ and $\nu$ are the same. 
If the tubular neighborhood of a fiber in a Seifert manifold
has parameter $\mu> 1$, then the fiber is a \emph{$\mu$-fold exceptional fiber}; 
if $\mu=1$, then the fiber is a \emph{regular fiber}.

\begin{thm} {\rm (Seifert~\cite{Seifert1933})}
Every Seifert manifold has at most finitely
many exceptional fibers.
\end{thm}

The set of all fibers in a Seifert manifold $M^3$, equipped with
the quotient topology, is a closed (orientable or nonorientable) 
$2$-dimensional surface, called the \emph{orbit surface} of $M^3$.

If $M^2$ is any surface, then $M^2\times S^1$ is a
Seifert manifold with fibers $\{x\}\times S^1$, $x\in M^2$, and orbit
surface $M^2$. Here, obviously, the orbit surface $M^2$ can be embedded
into $M^2\times S^1$, but this need not be the case in general.

If $M^3$ has $r$ exceptional fibers, then we can \emph{drill out}
these fibers by removing disjoint open tubular neighborhoods of the $r$ fibers.
This way we obtain a fibered space with $r$ toroidal boundary components.
We can close this space by \emph{filling in} regular solid tori
in order to obtain a Seifert manifold with no exceptional fibers that
still has the same orbit surface. However, there are, in general, 
different ways to do this, which requires to first classify 
all Seifert manifolds without exceptional fibers.

\begin{thm} {\rm (Seifert~\cite{Seifert1933})}
If from a Seifert manifold without exceptional fibers
one removes an open tubular neighborhood of an arbitrary
regular fiber, then the resulting ``classifying space''
(with toroidal boundary) is of one of six distinct possible types.
\end{thm}
The classifying spaces are denoted as $\{ Oo,g\}$, $\{ On,g\}$,
$\{ No,g\}$, $\{ NnI,g\}$, $\{ NnI\!I,g\}$, $\{ NnI\!I\!I,g\}$, where the capital
letters $O$ and $N$ stand for the orientability respectively nonorientability
of the total space and the small letters $o$ and $n$ for an 
orientable respectively nonorientable (punctured) orbit
surface of genus $g$. Note that, in general, there are three different
types $I$, $I\!I$, and $I\!I\!I$
of nonorientable classifying spaces in the case of a nonorientable orbit
surface of genus~$g$.

\bigskip

\emph{The type $\{ Oo,g\}$}. 
Orientable Seifert manifolds $M^3$ with orientable orbit surface have
a particularly simple description. Those without exceptional fibers
can be obtained as follows. Start with a (once) \emph{punctured} 
orientable surface of genus $g$, i.e., remove an open disk from the surface. 
Then form the topological direct product with a circle $S^1$.
The resulting fibered space is the classifying space $\{ Oo,g\}$ 
and has one toroidal boundary component.

Let $M$ be a simple closed oriented curve on the boundary torus $T^2$ of a solid torus $V$. 
If $M$ is nullhomotopic in $V$ but not contractible on $T^2$, then $M$
is called a \emph{meridian}. If, in addition, the solid torus $V$ is fibered, then
every simple closed curve on the boundary torus $T^2$ that intersects
each fiber exactly once is a \emph{crossing curve}.

In a next step we close $\{ Oo,g\}$ by gluing in a regular solid torus.
For this let $Q_0$ be the boundary circle of the (embedded) punctured
surface and let $H_0$
be a fiber of the fibered boundary torus of $\{ Oo,g\}$.
Furthermore, let $M_0$ be a simple closed curve on the boundary
torus homotopy equivalent to $Q_0+bH_0$, i.e.,
$$M_0\sim Q_0+bH_0,$$
where $b$ is any integer. Then there is a unique regular fibered solid torus $V$
whose boundary torus can be mapped under a fiber preserving map
onto the boundary torus of $\{ Oo,g\}$ such that $M_0$ is nullhomotopic
in $V$. 

\begin{figure}
\begin{center}
\footnotesize
\psfrag{H}{$H_0$}
\psfrag{M}{$M_0$}
\psfrag{Q}{$Q_0$}
\includegraphics[width=.265\linewidth]{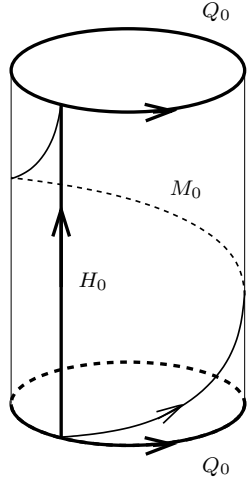}
\end{center}
\caption{Boundary torus with $M_0\sim Q_0+1\cdot H_0$.}
\label{fig:cylinder2}
\end{figure}

\begin{figure}
\begin{center}
\footnotesize
\psfrag{H}{$H_0$}
\psfrag{M}{$M_0$}
\psfrag{Q}{$Q_0$}
\includegraphics[width=.85\linewidth]{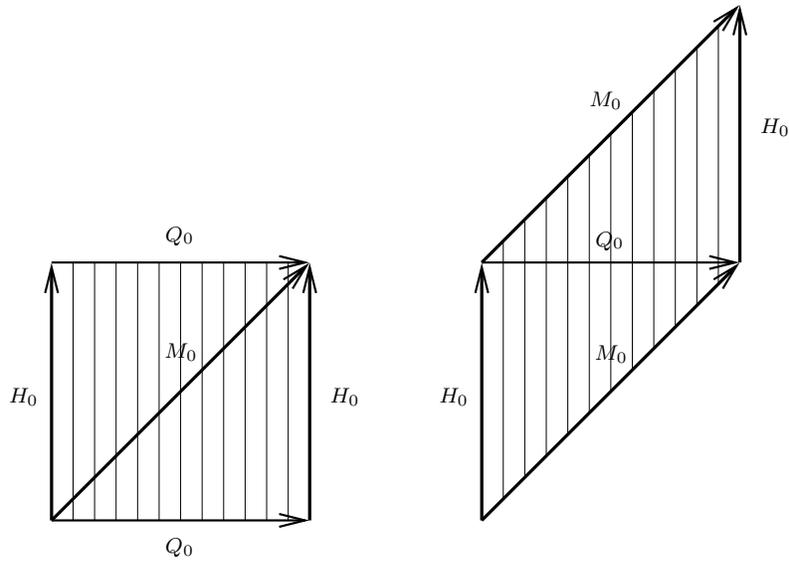}
\end{center}
\caption{Fibered boundary torus with $M_0\sim Q_0+1\cdot H_0$.}
\label{fig:cylinder3}
\end{figure}

In Figure~\ref{fig:cylinder2} and Figure~\ref{fig:cylinder3} (left) 
we display the boundary torus of $\{ Oo,g\}$
together with a curve $M_0\sim Q_0+1\cdot H_0$ that becomes a meridian
on the boundary of the corresponding fibered solid torus in Figure~\ref{fig:cylinder3} (right).

By closing the boundary torus of the classifying space $\{ Oo,g\}$ 
we obtain for every integer $b$ a different 
\emph{class space} $\{ Oo,g\mid b\}$ (without exceptional fibers).

Orientable Seifert manifolds with an orientable orbit surface of genus $g$
and with $r>0$ exceptional fibers 
can be constructed in a similar manner.
We start again with the orientable surface of genus $g$.
We remove $r+1$ disjoint open discs such that the resulting
surface has $r+1$ boundary circles $Q_0,Q_1,\ldots,Q_r$. Then we form 
the product with $S^1$, and get a $3$-manifold
with $r+1$ toroidal boundary components.
As above, we seal the boundary torus that has the crossing curve $Q_0$ 
with a regular torus corresponding to the parameter $b$. 
On the other boundary tori we consider curves 
$$M_i\sim \alpha_i Q_i+\beta_iH_i,\quad i=1,\ldots,r,$$
for fibers $H_1,\ldots,H_r$
and pairwise prime integers $\alpha_i$ and $\beta_i$, 
which, without loss of generality, can be chosen such that  
$$\alpha_i >1\hspace{5mm}\mbox{\textrm{and}}\hspace{5mm}0<\beta_i<\alpha_i.$$
Then for each such $M_i$ there is a unique fibered solid torus (with
parameters $\alpha_i$ and $\beta_i$) which has the curve $M_i$ as meridian.
By gluing in these tori we obtain a Seifert manifold
$\{ Oo,g\mid b;(\alpha_1,\beta_1),\ldots,(\alpha_r,\beta_r)\}$
with $r$ exceptional fibers determined by the parameters
$(\alpha_i,\beta_i)$, $i=1,\ldots,r$, and every Seifert manifold 
with classifying space $\{ Oo,g\}$ arises this way.

\bigskip

\emph{The type $\{ No,g\}$}. In this case the classifying space is
not the topological product of the punctured orientable surface
of genus $g$ with $S^1$ anymore. Rather, the orientable surface of
genus $g$ is cut into the corresponding fundamental polygon (see
\cite[Ch.\ 1.3]{Stillwell1980})
and punctured by cutting off the vertices.

\begin{figure}
\begin{center}
\footnotesize
\psfrag{a}{$a$}
\psfrag{b}{$b$}
\includegraphics[width=.5\linewidth]{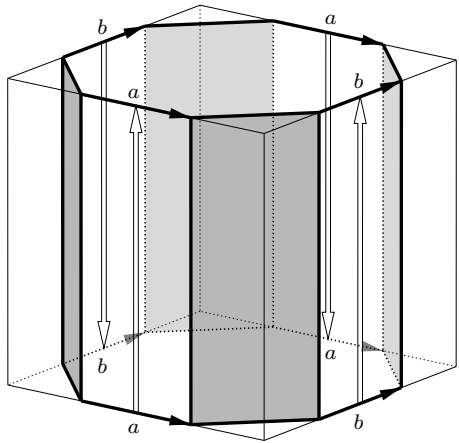}
\end{center}
\caption{Classifying space $\{ No,1\}$.}
\label{fig:SeifertNo1}
\end{figure}

The fundamental polygon of the orientable surface of genus 1, i.e.,
of the $2$-torus $T^2$, is a square with identified opposite oriented edges $a$ and $b$,
and is depicted as the top side (as well as the bottom side) of
the cube in Figure~\ref{fig:SeifertNo1}.

In fact, before corresponding edges of the punctured
fundamental polygon are identified, we take the product
of the polygon (with cut off vertices) with an interval $I$. 
In Figure~\ref{fig:SeifertNo1} this gives a vertical rectangular face 
for every edge of the original square and a shaded rectangular face 
for every cut off vertex of the square. Then the polygon on the top
side is identified with the copy on the bottom side
(so that we have the product of the polygon with cut off vertices with $S^1$).
In a second step, white rectangular side faces are pairwise identified
such that edge labels and arrows match. By flipping the side faces,
the resulting total space $\{ No,1\}$ becomes nonorientable.
Its boundary consists of the shaded faces that form a torus when glued together.

The same procedure can be carried out for the punctured
fundamental polygon of any orientable surface of genus $g$
and yields a nonorientable classifying space $\{ No,g\}$
with a toroidal boundary.

\bigskip

\begin{figure}
\begin{center}
\footnotesize
\psfrag{a}{$a$}
\psfrag{b}{$b$}
\psfrag{c}{$c$}
\includegraphics[width=.575\linewidth]{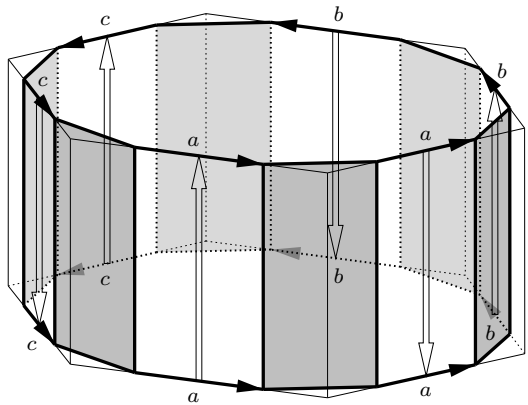}
\end{center}
\caption{Classifying space $\{ On,3\}$.}
\label{fig:SeifertOn3}
\end{figure}

\emph{The type $\{ On,g\}$}. As in the case $\{ No,g\}$ we begin
with the punctured fundamental polygon of the nonorientable
surface of genus $g$, form the product with $S^1$ and identify
side faces corresponding to edges of the fundamental polygon.
In order to obtain an orientable total space, all side faces need to
be flipped before identification. (See Figure~\ref{fig:SeifertOn3} for the
classifying space $\{ On,3\}$.)

\bigskip

\emph{The type $\{ NnI,g\}$}. If none of the side faces is flipped before identification,
then the resulting classifying space $\{ NnI,g\}$ is the product of the
punctured nonorientable surface of genus $g$ with $S^1$.

\bigskip

\begin{figure}
\begin{center}
\footnotesize
\psfrag{a}{$a$}
\psfrag{b}{$b$}
\psfrag{c}{$c$}
\includegraphics[width=.575\linewidth]{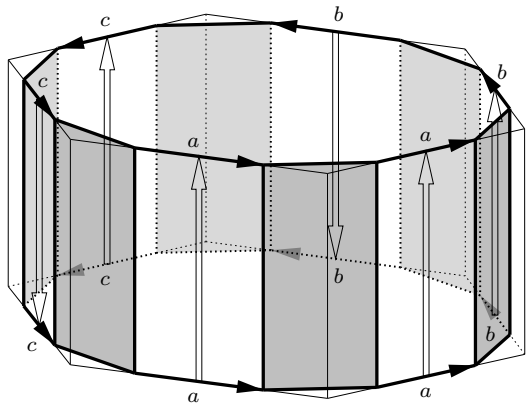}
\end{center}
\caption{Classifying space $\{ NnI\!I,3\}$.}
\label{fig:SeifertNnII3}
\end{figure}

\emph{The type $\{ NnI\!I,g\}$}. In this case, all side faces are
flipped before identification, except for one pair. This requires that
$g\geq 2$. (See Figure~\ref{fig:SeifertNnII3} for the
classifying space $\{ NnI\!I,3\}$.)

\bigskip

\emph{The type $\{ NnI\!I\!I,g\}$}. For a nonorientable surface of genus
$g\geq 3$ and at least two pairs but not all pairs of side faces
flipped before identification the resulting classifying space $\{
NnI\!I\!I,g\}$ is nonorientable but different from $\{ NnI\!I,g\}$.

\bigskip

In order to insert $r>0$ exceptional fibers into one of these classifying spaces
we proceed as in the case $\{ Oo,g\}$ and remove $r$ disjoint open discs 
from the starting polygon with cut off vertices before the product 
with $I$ and the identifications of the side faces are carried out. 
We then get a space with $r+1$ toroidal boundary
components. 

For the orientable space $\{ On,g\}$ the closing
procedure is as for $\{ Oo,g\}$. 
For the nonorientable 
spaces the parameters $\alpha_i$ and $\beta_i$ 
can be chosen (without loss of generality) such that  
$$\alpha_i >1\hspace{5mm}\mbox{\textrm{and}}\hspace{5mm}0<\beta_i\leq\frac{1}{2}\alpha_i.$$
If $\alpha_i >2$ for all $i=1,\ldots,r$, 
then the parameter $b$ can be restricted to $b=0$ or $b=1$. 
If there is at least one $i$ with $\alpha_i=2$, then only the case
$b=0$ has to be considered.

\pagebreak

\begin{thm} {\rm \cite{Seifert1933}}
Every Seifert manifold is, up to fiber (and orientation) preserving homeomorphisms, 
of one of the following uniquely determined types:\\[-1mm]

\begin{tabular}{@{\hspace{.2cm}}l}
$\{ Oo,g\mid b;(\alpha_1,\beta_1),\ldots,(\alpha_r,\beta_r)\}$, \\[1mm]
$\{ On,g\mid b;(\alpha_1,\beta_1),\ldots,(\alpha_r,\beta_r)\}$, 
\end{tabular}\\[3mm]
with\, $\alpha_i >1$,\, $0<\beta_i<\alpha_i$, and integer $b$, or \\

\begin{tabular}{@{\hspace{.2cm}}l}
$\{ No,g\mid b;(\alpha_1,\beta_1),\ldots,(\alpha_r,\beta_r)\}$,  \\[1mm]
$\{ NnI,g\mid b;(\alpha_1,\beta_1),\ldots,(\alpha_r,\beta_r)\}$, \\[1mm]
$\{ NnI\!I,g\mid b;(\alpha_1,\beta_1),\ldots,(\alpha_r,\beta_r)\}$,\hspace{3.6mm} $g\geq2$,\\[1mm]
$\{ NnI\!I\!I,g\mid b;(\alpha_1,\beta_1),\ldots,(\alpha_r,\beta_r)\}$,\quad $g\geq3$,\\
\end{tabular}\\[3mm]
with\, $\alpha_i >1$,\, $0<\beta_i\leq\frac{1}{2}\alpha_i$,\,
and\, $b=0$\, or\, $b=1$ (the latter only when all\, $\alpha_i >2$).
\end{thm}

\noindent
Reversing the orientation yields\\[-1mm]

\begin{tabular}{@{\hspace{.2cm}}l}
$\{ Oo,g\mid -r-b;(\alpha_1,\alpha_1-\beta_1),\ldots,(\alpha_r,\alpha_r-\beta_r)\}$\\[1mm]
\end{tabular}\\[3mm]
respectively\\[-1mm]

\begin{tabular}{@{\hspace{.2cm}}l}
$\{ On,g\mid -r-b;(\alpha_1,\alpha_1-\beta_1),\ldots,(\alpha_r,\alpha_r-\beta_r)\}$
\end{tabular}\\[3mm]
in the cases of orientable spaces.

It may occur that homeomorphic Seifert manifolds are defined by different sets
of Seifert invariants and cannot be mapped onto each other
by simply reversing the orientation (or by permuting the exceptional
fibers). For example, for $\alpha\geq 2$ the Seifert manifolds $\{ Oo,0\mid
0;(\alpha,1)\}$ 
are all homeomorphic to the sphere $S^3$.

Seifert manifolds contained in the following list
are called \emph{small} and were discussed by 
Orlik and Raymond \cite{OrlikRaymond1969} (see also \cite[Ch.\ 5.3]{Orlik1972}):\\[-1mm]

\begin{tabular}{@{\hspace{.2cm}}l}
$\{ Oo,0\mid b\}$,\, $\{ Oo,0\mid b;(\alpha_1,\beta_1)\}$,\, $\{ Oo,0\mid b;(\alpha_1,\beta_1),(\alpha_2,\beta_2)\}$,\\[1mm]
$\{ Oo,0\mid b;(\alpha_1,\beta_1),(\alpha_2,\beta_2),(\alpha_3,\beta_3)\}$\, with\, $\frac{1}{\alpha_1}+\frac{1}{\alpha_2}+\frac{1}{\alpha_3}>1$,\\[1mm]
$\{ Oo,0\mid -2;(2,1),(2,1),(2,1),(2,1)\}$,\\[1mm]
$\{ Oo,1\mid b\}$,\, $\{ No,1\mid b\}$,\\[1mm]
$\{ On,1\mid b\}$,\, $\{ NnI,1\mid b\}$,\, $\{ On,1\mid b;(\alpha_1,\beta_1)\}$,\, $\{ NnI,1\mid b;(\alpha_1,\beta_1)\}$,\\[1mm]
$\{ On,2\mid b\}$,\, $\{ NnI,2\mid b\}$,\, $\{ NnI\!I,2\mid b\}$.
\end{tabular}\\[1mm]

\noindent
Some of the small Seifert manifolds are homeomorphic 
although they have different sets of defining Seifert invariants. 
All Seifert manifolds not in the list are called \emph{large}.

Large Seifert manifolds were classified up to homeomorphism
by Orlik, Vogt, and Zieschang \cite{OrlikVogtZieschang1967}:

\begin{thm}{\rm \cite{OrlikVogtZieschang1967}}
If $M$ and $N$ are two large Seifert manifolds,
then the following statements are equivalent:
\begin{itemize}
\item[(a)] $M$ and $N$ have the same Seifert invariants (up to
  orientation reversion),
\item[(b)] $M$ and $N$ are homeomorphic,
\item[(c)] $M$ and $N$ have isomorphic fundamental groups.
\end{itemize}
\end{thm}
Moreover, it follows from work of Waldhausen (\cite{Waldhausen1967a}, \cite{Waldhausen1967b})
that every large Seifert manifold $M$ is \emph{irreducible}, i.e., every
embedded $2$-sphere in $M$ bounds a $3$-ball in $M$. 
In particular, all large Seifert manifolds are prime. 
With the exception of
$S^2\hbox{$\times\hspace{-1.62ex}\_\hspace{-.4ex}\_\hspace{.7ex}$}S^1$,
$S^2\!\times\!S^1$, and ${\mathbb R}{\bf P}^{\,3}\#\,{\mathbb R}{\bf P}^{\,3}$
also all small Seifert manifolds are irreducible.
Both $2$-sphere bundles over $S^1$ are prime, which leaves
${\mathbb R}{\bf P}^{\,3}\#\,{\mathbb R}{\bf P}^{\,3}$
as the only non-prime Seifert manifold.

\bigskip

Seifert manifolds appear in various settings, for example 
as complex hypersurfaces (see Section~\ref{sec:hypersurfaces})
or as quotients of the geometries $S^2\times{\mathbb R}^1$, $S^3$, $E^3$, ${\rm Nil}$,
$H^2\times{\mathbb R}^1$, and $\widetilde{SL}(2,{\mathbb R})$.
However, it took until the early 1980s
that a clear picture of the connections between these different concepts was reached.
The main link is given by two invariants for a Seifert manifold $M$,
the \emph{Euler number} $e(M)$, introduced by Neumann and Raymond \cite{NeumannRaymond1978},
and the \emph{Euler characteristic} $\chi (X)$ of the base orbifold $X$ of
Thurston \cite[Ch.~13]{Thurston2002}. If $M$ is orientable, then 
$$e(M):=-\big(b+\sum_{i=1}^r\frac{\beta_i}{\alpha_i}\big),$$
and\, $e(M)\!:=\!0$\, if $M$ is nonorientable. The Euler
characteristic of the base orbifold is given by
$$\chi(X):=\chi(S)-\sum_{i=1}^r\big(1-\frac{1}{\alpha_i}\big),$$
where $\chi(S)$ denotes the Euler characteristic of the orbit surface.
\begin{thm}{\rm (Scott \cite[5.3]{Scott1983})}
The geometric type of a Seifert manifold
depends only on the sign of the Euler
characteristic $\chi$ of the base orbifold and on
whether $e=0$ or $e\neq 0$.
\end{thm}
The corresponding geometry of a Seifert manifold can be read off 
from Table~\ref{tbl:seifert_geometries}. In particular, 
it follows that all $S^3$-, ${\rm Nil}$-, and $\widetilde{SL}(2,{\mathbb R})$-spaces
are orientable. 

\begin{table}
\small\centering
\defaultaddspace=0.3em
\caption{The six geometries of Seifert manifolds.}\label{tbl:seifert_geometries}
\begin{tabular}{@{}l@{\hspace{1.2cm}}l@{\hspace{1.2cm}}l@{\hspace{1.2cm}}l@{}}
\\
\toprule
 \addlinespace
                & $\chi>0$                 & $\chi=0$       & $\chi<0$ \\
 \addlinespace
\midrule
 \addlinespace
 \addlinespace
 $e=0$          & $S^2\times{\mathbb R}^1$ & $E^3$          & $H^2\times{\mathbb R}^1$ \\
 \addlinespace
 $e\neq 0$      & $S^3$                    & ${\rm Nil}$    & $\widetilde{SL}(2,{\mathbb R})$ \\
 \addlinespace
\bottomrule
\end{tabular}
\end{table}

\noindent

\bigskip

We will next discuss the different geometries, and we will construct triangulations
for a number of spaces. There are only four manifolds modeled on the
geometry $S^2\times{\mathbb R}^1$ and only ten manifolds
modeled on $E^3$. We will give triangulations for all of these
as well as for products and some twisted products of surfaces with $S^1$.
Triangulations of particular spherical $3$-manifolds will be obtained 
geometrically via pairwise identifying boundary faces of well-known polyhedra.

Alternative to these constructions, an algorithmic procedure that yields
concrete triangulations for \emph{all} Seifert manifolds will
be given in \cite{BrehmLutz2002pre}. These triangulations
can easily be produced with the GAP-program SEIFERT \cite{SEIF},
which requires as input only the Seifert invariants.


\mathversion{bold}
\section{($S^2\times{\mathbb R}$)-Spaces}
\mathversion{normal}

Seifert \cite{Seifert1931} described in 1931 four quotient spaces 
of $S^2\times{\mathbb R}$. Tollefson \cite{Tollefson1974} then
showed in 1974 that these are already all ($S^2\times{\mathbb R}$)-spaces.

\begin{thm} {\rm (Tollefson \cite{Tollefson1974})}
There are precisely four $3$-manifolds modeled on $S^2\!\times\!{\mathbb R}$,
namely\, $S^2\!\times\!S^1$, $S^2\hbox{$\times\hspace{-1.62ex}\_\hspace{-.4ex}\_\hspace{.7ex}$}S^1$, 
${\mathbb R}{\bf P}^{\,2}\!\times S^1$, and ${\mathbb R}{\bf P}^{\,3}\#\,{\mathbb R}{\bf P}^{\,3}$.
\end{thm}
 
In this section, we will see how triangulations for products 
and twisted products can be obtained. We will also list the Seifert
invariants for the four ($S^2\times{\mathbb R}$)-spaces.

\bigskip

\emph{Products}. Let $M$ be a triangulable $m$-dimensional and 
$N$ be a triangulable $n$-dimensional manifold. Then their topological
product $M\times N$ can be triangulated as follows.
We start with a triangulation of $M$ and a triangulation of $N$.
In a first step, we form for every $m$-dimensional facet $\Delta_m$ of $M$ 
the direct product with each $n$-dimensional facet $\Delta_n$ of~$N$.
The union of all these products of simplices $\Delta_m\times \Delta_n$
gives a cell-decomposition of $M\times N$
which then will be triangulated consistently in a second step.

\begin{figure}
\begin{center}
\footnotesize
\psfrag{(v_0,w_0)}{$(v_0,w_0)$}
\psfrag{(v_n,w_0)}{$(v_n,w_0)$}
\psfrag{(v_0,w_n)}{$(v_0,w_n)$}
\psfrag{(v_n,w_n)}{$(v_n,w_n)$}
\includegraphics[width=.55\linewidth]{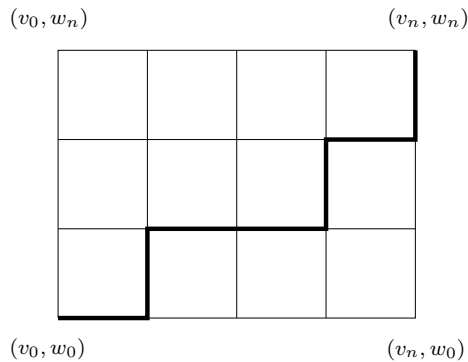}
\end{center}
\caption{Staircase triangulation of $\Delta_m\times \Delta_n$.}
\label{fig:staircase}
\end{figure}
One particular triangulation of the direct product
$\Delta_m\times \Delta_n$ of an $m$-simplex $\Delta_m$
having vertices $v_0,\ldots,v_m$ with an $n$-simplex $\Delta_n$
having vertices $w_0,\ldots,w_n$ is the \emph{staircase triangulation}
described in \cite{BilleraCushmanSanders1988},
\cite[Ch.\ 7]{GelfandKapranovZelevinsky1994},
\cite{Lee1997}, and \cite{Santos2000}. (See also \cite[Ch.\ II \S 8]{EilenbergSteenrod1952}.)
The product $\Delta_m\times \Delta_n$ has 
the vertex set 
$$\{ (v_k,w_l) \mid 0\leq k\leq m,\, 0\leq l\leq n\},$$
which we identify with the vertex set of a rectangular grid, as in Figure~\ref{fig:staircase}.
The $m+n+1$ vertices of every (monotone increasing) lattice path from $(v_0,w_0)$ to $(v_m,w_n)$
define an $(m+n)$-dimensional facet of $\Delta_m\times \Delta_n$.
For example, there are altogether three lattice paths for the product $\Delta_2\times \Delta_1$,
which give a triangulation of the prism $\Delta_2\times \Delta_1$ into
the following three tetrahedra, $\{(v_0,w_0),(v_1,w_0),(v_2,w_0),(v_2,w_1)\}$, 
$\{(v_0,w_0),(v_1,w_0),(v_1,w_1),(v_2,w_1)\}$, 
and $\{(v_0,w_0),(v_0,w_1),(v_1,w_1),(v_2,w_1)\}$.

If the vertex sets of the triangulations of $M$ respectively $N$ 
are totally ordered, then the union of the facets of the staircase triangulations
of all the cells $\Delta_m\times \Delta_n$ of $M\times N$
yields a consistent \emph{product triangulation} of $M\times N$.
If the triangulations of $M$ and $N$ have $s$ and $t$ vertices and
$u$ and $v$ facets, respectively, then the product triangulation of $M\times N$
has $s\cdot t$ vertices and $u\cdot v\cdot (m+n)!/(m!n!)$ maximal
faces \cite{Lee1997}.

\begin{figure}
\begin{center}
\footnotesize
\psfrag{1}{1}
\psfrag{2}{2}
\psfrag{3}{3}
\psfrag{4}{4}
\psfrag{5}{5}
\psfrag{6}{6}
\psfrag{7}{7}
\psfrag{8}{8}
\psfrag{9}{9}
\psfrag{10}{10}
\psfrag{11}{11}
\psfrag{12}{12}
\includegraphics[width=.675\linewidth]{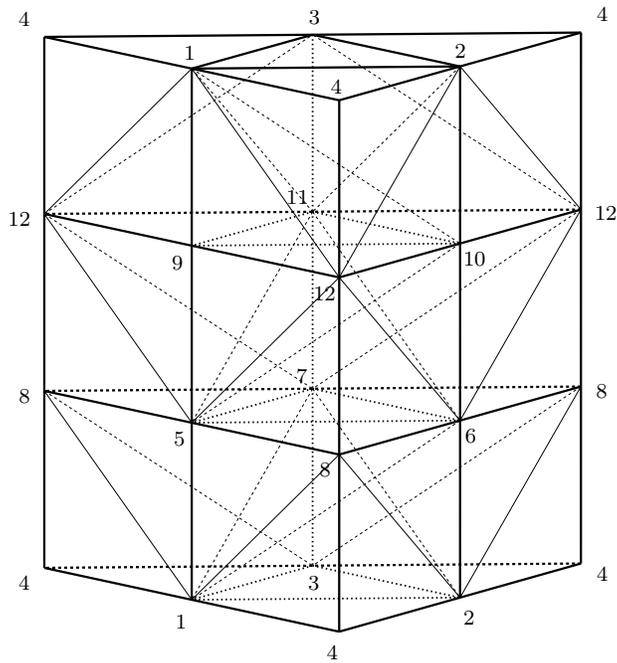}
\end{center}
\caption{Product triangulation of $S^2\!\times\!S^1$.}
\label{fig:product_triangulation}
\end{figure}

In Figure~\ref{fig:product_triangulation} we depict the product
triangulation of $S^2\!\times\!S^1$ where $S^2$ is triangulated 
as the boundary of the tetrahedron $1234$. The four triangles $123$,
$124$, $134$, and $234$ can be seen at the top and bottom of the figure.
The circle $S^1$ is triangulated with three edges. 
For simplicity, we have not labeled the vertices of the product
triangulation of $S^2\!\times\!S^1$ in the figure with pairs of vertices of $S^2$ and $S^1$.
Instead, we used a set of $4$ new vertices for every new ``horizontal
$S^2$-layer''. With respect to this labeling the rectangular faces of
the prisms are triangulated using the respective diagonal which contains the vertex with
the smallest label in the rectangle. The triangulation of the three rectangular faces of
every prism induces a unique triangulation of the prism itself into
three tetrahedra.

\bigskip

\emph{Twisted Products over $S^1$}. Once we know how to
triangulate direct products, then twisted products over $S^1$
are easy to describe.

Let $M$ be a triangulated $d$-dimensional manifold with $m$ vertices
and let $I$ be an interval triangulated with four vertices. The
product triangulation of $M\!\times I$ has $4\cdot m$ vertices
and has as its boundary two disjoint copies of $M$. If we glue these
copies of $M$ together trivially by identifying the respective
vertices, then we obtain a triangulation of $M\!\times\! S^1$ 
with $3\cdot m$ vertices. 

\begin{figure}
\begin{center}
\footnotesize
\psfrag{1}{1}
\psfrag{2}{2}
\psfrag{3}{3}
\psfrag{4}{4}
\psfrag{5}{5}
\psfrag{6}{6}
\psfrag{7}{7}
\psfrag{8}{8}
\psfrag{9}{9}
\psfrag{10}{10}
\psfrag{11}{11}
\psfrag{12}{12}
\psfrag{13=2}{13=2}
\psfrag{14=1}{14=1}
\psfrag{15=3}{15=3}
\psfrag{16=4}{16=4}
\includegraphics[width=.705\linewidth]{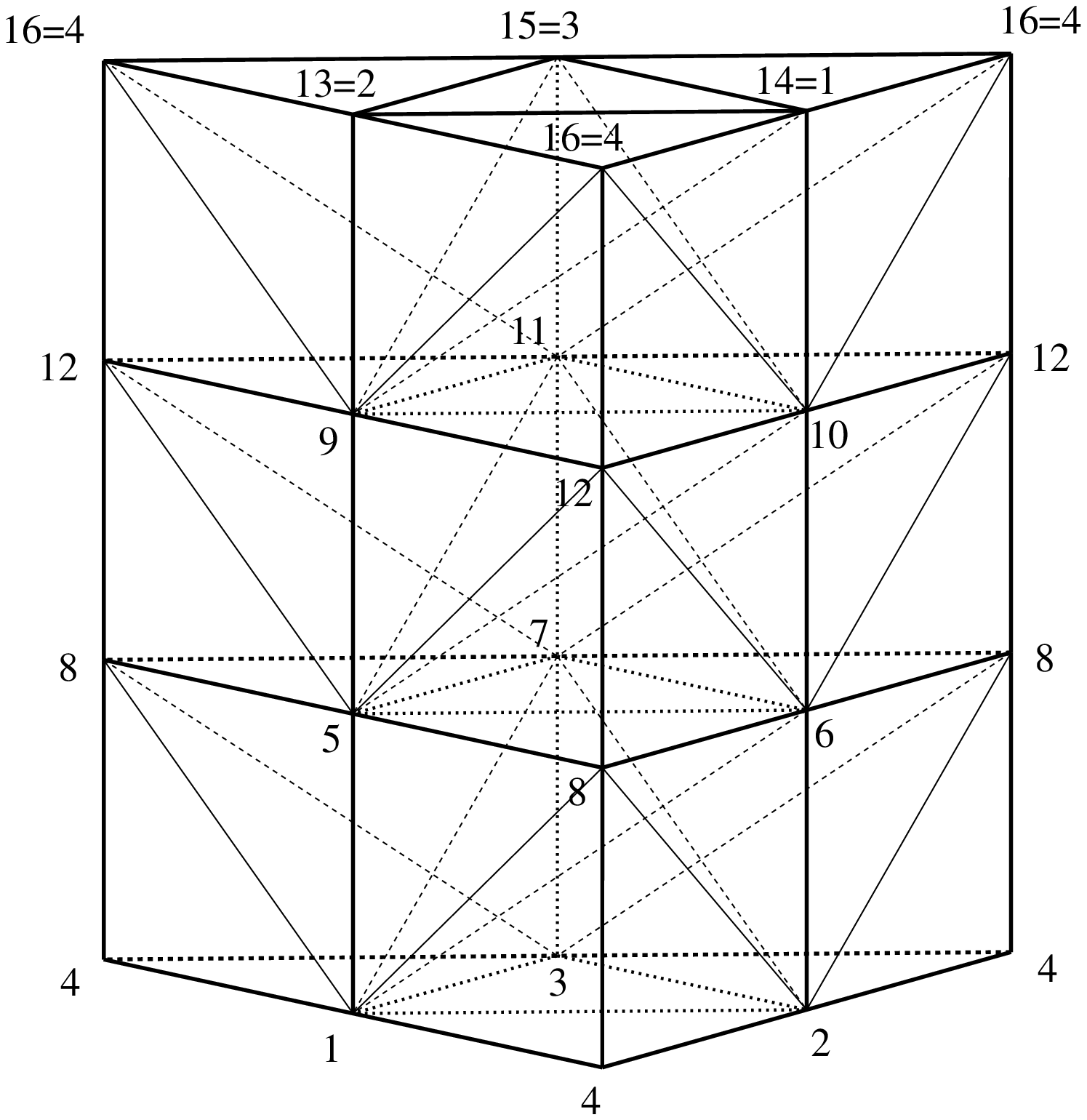}
\end{center}
\caption{Triangulation of the twisted product $S^2\hbox{$\times\hspace{-1.62ex}\_\hspace{-.4ex}\_\hspace{.7ex}$}S^1$.}
\label{fig:twist_triangulation}
\end{figure}

Instead of the trivial identification we can first apply any 
combinatorial automorphism of $M$ to one of the boundary components
of $M\!\times I$ before gluing the components together.
This way we obtain \emph{twisted product triangulations} 
of $M$ over $S^1$. Note that the homeomorphism type of the resulting
manifolds will, in general, depend on the combinatorial automorphism 
that is used for the gluing. Topologically, this construction
yields the \emph{mapping torus} of a homeomorphism.

In Figure~\ref{fig:twist_triangulation} the vertices 
$13$ and $14$ of the top copy of $S^2$ are flipped before the
identification with the bottom copy, so the two 
copies of $S^2$ are glued together under a map that reverses the orientation of $S^2$.
The resulting space is the nonorientable twisted product 
$S^2\hbox{$\times\hspace{-1.62ex}\_\hspace{-.4ex}\_\hspace{.7ex}$}S^1$.
There are exactly two $S^2$-bundles over $S^1$ \cite[\S 26]{Steenrod1951}.
Thus we only need to take care whether we change the orientation
of $S^2$ under the gluing or not.

If we start with any surface $M^2$ with $f$-vector $f=(f_0,f_1,f_2)$,
then the product triangulation of $M^2\!\times\!S^1$
or any twisted product triangulation of $M^2$ over $S^1$ has
$f$-vector $f=(3f_0,3f_0+6f_1,6f_1+4f_2,4f_2)$.

\bigskip

Let us next analyze which Seifert invariants 
can occur for ($S^2\times{\mathbb R}$)-spaces.
The Euler characteristic of the base orbifold $X$ has to be greater than $0$, i.e.,
$$\chi (X)=\chi(S) -\sum_{i=1}^r\big(1-\frac{1}{\alpha_i}\big)>0.$$
It follows that the orbit surface $S$ has to be
$S^2$ with $\chi(S^2)=2$, or ${\mathbb R}{\bf P}^{\,2}$ with
$\chi({\mathbb R}{\bf P}^{\,2})=1$, since all other surfaces have
Euler characteristic $\chi(S)\leq 0$ and since
$1>(1-\frac{1}{\alpha_i})\geq\frac{1}{2}$ for every exceptional fiber.
It also follows that the number of exceptional fibers $r$ is at most $3$
for the orbit surface $S^2$ and at most $1$ 
for ${\mathbb R}{\bf P}^{\,2}$.

If the orbit surface is $S^2$, then the total space $M$ is orientable.
In particular, we need that
$$e(M)=-\big(b+\sum_{i=1}^r\frac{\beta_i}{\alpha_i}\big)=0.$$
If $r=0$, then necessarily $b=0$.
The case $r=1$ cannot occur, since 
$\frac{\beta}{\alpha}$ is not an integer 
(thus $e(M)=-(b+\frac{\beta}{\alpha})\neq 0$).
For $r=2$ it is $0<\frac{\beta_1}{\alpha_1}+\frac{\beta_2}{\alpha_2}<2$,
so the sum $\frac{\beta_1}{\alpha_1}+\frac{\beta_2}{\alpha_2}$ 
has to be equal to $1$, which can only be achieved
for $\alpha_2=\alpha_1$ and $\beta_2=\alpha_1-\beta_1$.
In the case $r=3$ we have for the base orbifold $X$ that $\chi (X)>0$
is equivalent to
$$\frac{1}{\alpha_1}+\frac{1}{\alpha_2}+\frac{1}{\alpha_3}>1.$$
There are exactly four types of \emph{platonic triples}
of integers $\alpha_i$ for which the inequality holds: $(2,2,\alpha_3)$,
$(2,3,3)$, $(2,3,4)$, and $(2,3,5)$.
However, in the first case 
$\sum_{i=1}^3\frac{\beta_i}{\alpha_i}=\frac{1}{2}+\frac{1}{2}+\frac{\beta_3}{\alpha_3}$
is not an integer and this also can be deduced for the other platonic
triples. Therefore, $r=3$ is not possible.

Let the orbit surface be ${\mathbb R}{\bf P}^{\,2}$.
If the total space is orientable, then the case $r=1$ can be ruled out
as above. Hence, $r=0$ and $b=0$. 
For a nonorientable total space $e(M)=0$, so $r\leq 1$
is the only restriction.

The topological types of the resulting Seifert manifolds
were determined by Orlik and Raymond \cite{OrlikRaymond1969}.
In Table~\ref{tbl:small_S2xR_spaces} 
we give the four ($S^2\times{\mathbb R}$)-spaces
together with their Seifert invariants and the
$f$-vectors of their smallest known triangulations.
The triangulations of $S^2\hbox{$\times\hspace{-1.62ex}\_\hspace{-.4ex}\_\hspace{.7ex}$}S^1$
and $S^2\!\times\!S^1$ are vertex-minimal. They can, for example, be obtained
from the product and the twisted product triangulation of $S^2$ over $S^1$
by applying bistellar flips to these. 

\begin{thm}
The four ($S^2\times{\mathbb R}$)-spaces have small triangulations
with at most $15$ vertices.
\end{thm}

More details on small triangulations of $3$-manifolds with at most $15$ vertices
are given in \cite{Lutz1999}.

\begin{table}
\small\centering
\defaultaddspace=0.3em
\caption{Smallest known triangulations of ($S^2\times{\mathbb R}$)-spaces.}\label{tbl:small_S2xR_spaces}
\begin{tabular*}{\linewidth}{@{\extracolsep{\fill}}l@{\hspace{2mm}}l@{\hspace{2mm}}l@{\hspace{2mm}}l@{\hspace{2mm}}l@{}}
\\
\toprule
 \addlinespace
    Manifold         & Seifert Fibrations & Homology & $f$-Vector \\
\midrule
 \addlinespace
 \addlinespace
 $S^2\hbox{$\times\hspace{-1.62ex}\_\hspace{-.4ex}\_\hspace{.7ex}$}S^1$ & $\{ NnI,1\mid 1\}$  
                     & $({\mathbb Z},{\mathbb Z},{\mathbb Z}_2,0)$ & (9,36,54,27) \\
                     & $\{ NnI,1\mid b;(\alpha,\beta)\}$\, with\, $b\alpha+\beta$\, odd \\ 
 \addlinespace
 $S^2\!\times\!S^1$  & $\{ Oo,0\mid 0\}$, & $({\mathbb Z},{\mathbb Z},{\mathbb Z},{\mathbb Z})$ & (10,40,60,30) \\
                     & $\{ Oo,0\mid -1;(\alpha,\beta),(\alpha,\alpha-\beta)\}$ & & \\
 \addlinespace
 ${\mathbb R}{\bf P}^{\,2}\!\times S^1$ & $\{ NnI,1\mid 0\}$ & $({\mathbb Z},{\mathbb Z}\oplus{\mathbb Z}_2,{\mathbb Z}_2,0)$ & (14,84,140,70) \\
                     & $\{ NnI,1\mid b;(\alpha,\beta)\}$\, with\, $b\alpha+\beta$\, even \\ 
 \addlinespace
 ${\mathbb R}{\bf P}^{\,3}\#\,{\mathbb R}{\bf P}^{\,3}$ & $\{ On,1\mid 0\}$ & $({\mathbb Z},{\mathbb Z}_2^{\,2},0,{\mathbb Z})$ & (15,86,142,71) \\
 \addlinespace
\bottomrule
\end{tabular*}
\end{table}


\section{Spherical Spaces}

Spherical $3$-dimensional Clifford-Klein manifolds
were classified essentially by Hopf \cite{Hopf1926} in
1926 and were described in detail by Threlfall and Seifert 
in the early 1930s (\cite{ThrelfallSeifert1931}, \cite{ThrelfallSeifert1933}).
Before that time all the known examples
of spherical $3$-dimensional spaces were, in fact, quotients $S^3/G$ 
of the Lie group $S^3\cong SU(2)$ with respect to left (or right) 
multiplication of finite subgroups $G$ of $S^3$.
The group $S^3$ is the simply connected double cover of 
the three-dimensional rotation group $SO(3)\cong{\mathbb R}{\bf P}^{\,3}$.
Finite subgroups of $SO(3)$ were already determined by Klein (\cite{Klein1884}, \cite{Klein1993}). 
These are 
the cyclic groups ${\mathbb Z}_n$ of order $n$, 
the dihedral groups $D_n$ of order $2n$, 
the tetrahedral group $T$ of order $12$, 
the octahedral group $O$ of order $24$, and 
the icosahedral group $I$ of order~$60$.

The double cover of ${\mathbb Z}_n$ in $S^3$ is the cyclic group ${\mathbb Z}_{2n}$,
the double covers of the other groups are the respective binary
groups; see Table~\ref{tbl:groups_spherical}.

\begin{table}
\small\centering
\defaultaddspace=0.3em
\caption{Groups with fixed point free orthogonal action on $S^3$.}\label{tbl:groups_spherical}
\begin{tabularx}{\linewidth}{@{}X@{}}
\\
\toprule
 \addlinespace
 \addlinespace
 $D^*_n=\{ x,y\mid x^2=(xy)^2=y^n\}$, binary dihedral group of order $4n$ \\
 \addlinespace
 $T^*=\{ x,y\mid x^2=(xy)^3=y^3,\, x^4=1\}$, binary tetrahedral group of order 24 \\
 \addlinespace
 $O^*=\{ x,y\mid x^2=(xy)^3=y^4,\, x^4=1\}$, binary octahedral group of order 48 \\
 \addlinespace
 $I^*=\{ x,y\mid x^2=(xy)^3=y^5,\, x^4=1\}$, binary icosahedral group of order 120 \\
 \addlinespace
 $D'_{2^k(2n+1)}=\{ x,y\mid x^{2^k}=1,\, y^{2n+1}=1,\, xyx^{-1}=y^{-1}\}$, where\, $k\geq 2$, $n\geq 1$, \\[1mm]
 \multicolumn{1}{@{}r@{}}{of order\, $2^k(2n+1)$} \\
 \addlinespace
 $T'_{8\cdot 3^k}=\{ x,y,z\mid x^2=(xy)^2=y^2,\, zxz^{-1}=y,\, zyz^{-1}=xy,\, z^{3^k}=1\}$, where\, $k\geq 1$, \\[1mm]
 \multicolumn{1}{@{}r@{}}{of order\, $8\cdot 3^k$} \\
 \addlinespace
\bottomrule
\end{tabularx}
\end{table}

The full group of orientation-preserving isometries of $S^3$ is $SO(4)$.
(Since orientation-reversing isometries have fixed points, these are
not of interest here.) Discrete subgroups $G$ of $SO(4)$ 
that act freely on $S^3$ are finite. Spherical $3$-manifolds that 
arise as quotients $S^3/G$ therefore have finite fundamental group $G$,
and for their classification of spherical $3$-manifolds,
Hopf, Threlfall and Seifert determined all finite subgroups of $SO(4)$ 
that can act freely on $S^3$. (See also Milnor \cite{Milnor1957}.)

\begin{thm} {\rm (Hopf \cite{Hopf1926}, Threlfall and Seifert (\cite{ThrelfallSeifert1931}, \cite{ThrelfallSeifert1933}))}
\label{thm:HopfThrelfallSeifert}
Any finite subgroup of $SO(4)$ that has a free action on $S^3$ is
either a cyclic group ${\mathbb Z}_n$, one of the groups 
$D^*_n$, $T^*$, $O^*$, $I^*$, $D'_{2^k(2n+1)}$, or $T'_{8\cdot 3^k}$,
given in Table~\ref{tbl:groups_spherical},
or a direct product of any of these groups with a cyclic group of
relatively prime order.
\end{thm}

Spherical $3$-manifolds $M$ are Seifert manifolds which have
Euler number $e(M)\neq 0$ and base orbifold $X$ 
of positive Euler characteristic $\chi(X)>0$. In particular,
they are orientable with orientable orbit surface $S^2$ or with
nonorientable orbit surface ${\mathbb R}{\bf P}^{\,2}$.
According to the analysis of the previous section,
fibered spherical $3$-manifolds with orbit surface $S^2$ or ${\mathbb R}{\bf P}^{\,2}$ 
can have at most three respectively one exceptional fibers. A complete list
of the corresponding spaces with all their possible fiberings
is given in Table~\ref{tbl:ssm}.
The spaces themselves are discussed below.

\begin{landscape}

\small
\defaultaddspace=.3em

\setlength{\LTleft}{0pt}
\setlength{\LTright}{0pt}
\begin{longtable}{@{\extracolsep{\fill}}lll@{}}
\caption{\protect\parbox[t]{15cm}{Spherical Seifert manifolds.}}\label{tbl:ssm}
\\
\toprule
 \addlinespace
    Manifold             & Seifert Fibrations & Fundamental Group\\
\midrule
\endfirsthead
\caption{\protect\parbox[t]{15cm}{Spherical Seifert manifolds.}}
\\
\toprule
 \addlinespace
    Manifold             & Seifert Fibrations & Fundamental Group\\
\midrule
\endhead
\bottomrule
\endfoot
 \addlinespace
 \addlinespace
 $S^3$                   & $\{ Oo,0\mid 1\}$, $\{ Oo,0\mid -1\}$,  & 0\\
                         & $\{ Oo,0\mid 0;(\alpha,1)\}$, \\
                         & $\{ Oo,0\mid -1;(\alpha,\alpha-1)\}$, \\
                         & $\{ Oo,0\mid b;(\alpha_1,\beta_1),(\alpha_2,\beta_2)\}$ \\
                         & with $(\alpha_1,\alpha_2)=1$ \\
                         & and $|b\alpha_1\alpha_2+\alpha_1\beta_2+\alpha_2\beta_1|=1$ \\
 \addlinespace
 \addlinespace
 $L(b,1)$                & $\{ Oo,0\mid b\}$,\, $|b|>1$ & ${\mathbb Z}_{|b|}$\\
 \addlinespace
 \addlinespace
 $L(b\alpha+\beta,\alpha')$  & $\{ Oo,0\mid b;(\alpha,\beta)\}$ & ${\mathbb Z}_{|b\alpha+\beta|}$ \\
                         & with\, $|b\alpha+\beta|>1$\\
                         & and\, $\alpha'\equiv\alpha\mod(b\alpha+\beta)$\\
                         & with\, $0<\alpha'<b\alpha+\beta$ \\
 \addlinespace
 \addlinespace
 $L(b\alpha_1\alpha_2+\alpha_1\beta_2+\alpha_2\beta_1,m\alpha_2-n\beta_2)$  & $\{ Oo,0\mid b;(\alpha_1,\beta_1),(\alpha_2,\beta_2)\}$ & ${\mathbb Z}_{|b\alpha_1\alpha_2+\alpha_1\beta_2+\alpha_2\beta_1|}$\\
                         & with $|b\alpha_1\alpha_2+\alpha_1\beta_2+\alpha_2\beta_1|>1$ \\
                         & and \, $m\alpha_1-n(b\alpha_1+\beta_1)=1$ \\
 \addlinespace
 \addlinespace
 $P(r)$                  & $\{ Oo,0\mid -1;(2,1),(2,1),(r,1)\}$ & $D^*_{r}$ \\
 \addlinespace
 \addlinespace
 generalized prism space & $\{ Oo,0\mid b;(2,1),(2,1),(\alpha_3,\beta_3)\}$ & ${\mathbb Z}_{|(b+1)\alpha_3+\beta_3|}\times D^*_{\alpha_3}$\, if\, $((b+1)\alpha_3+\beta_3,2\alpha_3)=1$, \\
                         && ${\mathbb Z}_{|m|}\times D'_{2^{k+2}\alpha_3}$\, for\, $(b+1)\alpha_3+\beta_3=2^{k}m$\\
                         && with\, $k\geq 1$\, and\, $(m,2)=1$, \\
                         & $\{ OnI,1\mid b;(\alpha_1,\beta_1)\}$ & ${\mathbb Z}_{\alpha_1}\times D^*_{|b\alpha_1+\beta_1|}$\, if\, $\alpha_1$\, is odd, \\
                         & with\, $|b\alpha_1+\beta_1|\neq 0$ & ${\mathbb Z}_{\alpha_1'}\times D'_{2^{k+2}|b\alpha_1+\beta_1|}$\, for\, $\alpha_1=2^k\alpha_1'$ \\
                         && with\, $k\geq 1$\, and\, $(\alpha_1',2)=1$\\
 \addlinespace
 \addlinespace
 $S^3/T^*$               & $\{ Oo,0\mid -1;(2,1),(3,1),(3,1)\}$ & $T^*$\\
 \addlinespace
 \addlinespace
 generalized octahedral space & $\{ Oo,0\mid b;(2,1),(3,\beta_2),(3,\beta_3)\}$ & 
                           ${\mathbb Z}_{|6b+3+2\beta_2+2\beta_3|}\times T^*$\, if\, $(6b+3+2\beta_2+2\beta_3,12)=1$, \\
                         && ${\mathbb Z}_{|m|}\times T'_{8\cdot 3^{k+1}}$\, for\, $6b+3+2\beta_2+2\beta_3=3^km$\\
                         && with\, $k\geq 1$\, and\, $(m,12)=1$ \\
 \addlinespace
 \addlinespace
 $S^3/O^*$               & $\{ Oo,0\mid -1;(2,1),(3,1),(4,1)\}$ & $O^*$ \\
 \addlinespace
 \addlinespace
 generalized truncated cube space & $\{ Oo,0\mid b;(2,1),(3,\beta_2),(4,\beta_3)\}$ & ${\mathbb Z}_{|12b+6+4\beta_2+3\beta_3|}\times O^*$ \\
 \addlinespace
 \addlinespace
 $S^3/I^*=\Sigma (2,3,5)$  & $\{ Oo,0\mid -1;(2,1),(3,1),(5,1)\}$ & $I^*$ \\
 \addlinespace
 \addlinespace
 generalized spherical dodecahedral space & $\{ Oo,0\mid b;(2,1),(3,\beta_2),(5,\beta_3)\}$ & ${\mathbb Z}_{|30b+15+10\beta_2+6\beta_3|}\times I^*$ \\
 \addlinespace
 \addlinespace
\end{longtable}

\end{landscape}

\begin{figure}
\begin{center}
\footnotesize
\psfrag{1}{1}
\psfrag{2}{2}
\psfrag{p}{$p$}
\psfrag{p-1}{$p-1$}
\psfrag{p-q-1}{$p-q-1$}
\psfrag{p-q}{$p-q$}
\psfrag{p-q+1}{$p-q+1$}
\psfrag{p-q+2}{$p-q+2$}
\includegraphics[width=.7\linewidth]{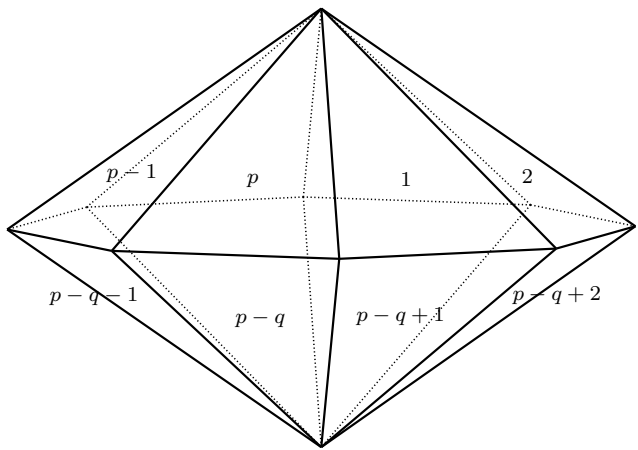}
\end{center}
\caption{The lens space $L(p,q)$.}
\label{fig:lens}
\end{figure}

Fibered spherical $3$-manifolds over $S^2$ with $r\leq 2$ exceptional fibers
are either the sphere $S^3$ or a lens space $L(p,q)$.
The \emph{lens space} $L(p,q)$ of Tietze \cite{Tietze1908} 
is defined as a solid lens with $p$ slices 
where top triangles are identified with bottom triangles under a twist
of $2\pi(q/p)$ as depicted in Figure~\ref{fig:lens}. Note that
$L(2,1)$ is the real projective space ${\mathbb R}{\bf P}^{\,3}$.

The homeomorphism classification of lens spaces
is subtle. A lens space $L(p,q)$ has fundamental group ${\mathbb Z}_{p}$.
However, two lens spaces $L(p,q)$ and $L(p,q')$ with different $q$ and
$q'$ need not be homeomorphic, not even homotopy equivalent.
\begin{thm} {\rm (Whitehead \cite{Whitehead1941})}
$L(p,q)\simeq L(p,q')$ if and only if\, $qq'\equiv\pm n^2\pmod p$ for
some integer $n$.
\end{thm}
For example, $L(5,1)$ and $L(5,2)$ are not homotopy equivalent,
although they have the same fundamental group (and the same homology
groups).
\begin{thm} {\rm (Reidemeister \cite{Reidemeister1936}, Brody \cite{Brody1960})}
$L(p,q)\cong L(p,q')$ if and only if either\, $q'\equiv\pm q\pmod p$\, or\,
$qq'\equiv\pm 1\pmod p$.
\end{thm}
As a consequence, for example, the lens spaces $L(7,1)$ and $L(7,2)$ are homotopy equivalent but 
not homeo\-mor\-phic.

Apart from lens spaces, two spherical $3$-manifolds
are homeomorphic if and only if they have the
same fundamental group.

Triangulations of lens spaces $L(p,q)$ with few simplices were constructed by Brehm and
\'Swiatkowski \cite{BrehmSwiatkowski1993}. 
They also gave a series of $D_{2(p+2)}$-sym\-metric triangulations
$S_{2(p+2)}$\, of\, $L(p,1)$\, on\, $2p+7$\, vertices. 
We applied bistellar flips to some of these examples. 
The $f$-vectors of the resulting small triangulations
are given in Table~\ref{tbl:small_lens_spaces}. For small
triangulations of other lens spaces see~\cite{BrehmLutz2002pre}.

\begin{table}
\small\centering
\defaultaddspace=0.3em
\caption{Smallest known triangulations of lens spaces $L(p,1)$.}\label{tbl:small_lens_spaces}
\begin{tabular}{@{}l@{\hspace{10mm}}l@{}}
\\
\toprule
 \addlinespace
    Manifold         & $f$-Vector \\ 
\midrule
 \addlinespace
 \addlinespace
 ${\mathbb R}{\bf P}^{\,3}=L(2,1)$
                     & (11,51,80,40) \\
 \addlinespace
 $L(3,1)$            & (12,66,108,54) \\
 \addlinespace
 $L(4,1)$            & (14,84,140,70) \\
 \addlinespace
 $L(5,1)$            & (15,97,164,82) \\
 \addlinespace
 $L(6,1)$            & (16,110,188,94) \\
 \addlinespace
 $L(7,1)$            & (17,123,212,106) \\
 \addlinespace
 $L(8,1)$            & (18,135,234,117) \\
 \addlinespace
 $L(9,1)$            & (18,144,252,126) \\
 \addlinespace
 $L(10,1)$           & (19,156,274,137) \\
 \addlinespace
\bottomrule
\end{tabular}
\end{table}

In the case of fibered spherical $3$-manifolds with three exceptional fibers 
$(\alpha_1,\beta_1)$, $(\alpha_2,\beta_2)$, and $(\alpha_3,\beta_3)$  
it follows from the preceding section that\linebreak
$\frac{1}{\alpha_1}+\frac{1}{\alpha_2}+\frac{1}{\alpha_3}>1$,
which is only possible for the platonic triples
$(2,2,\alpha_3)$, $(2,3,3)$, $(2,3,4)$, and $(2,3,5)$.

\bigskip

\begin{figure}
\begin{center}
\footnotesize
\psfrag{1}{1}
\psfrag{2}{2}
\psfrag{3}{3}
\psfrag{4}{4}
\psfrag{5}{5}
\psfrag{6}{6}
\psfrag{7}{7}
\psfrag{8}{8}
\psfrag{9}{9}
\psfrag{10}{10}
\psfrag{11}{11}
\psfrag{12}{12}
\psfrag{13}{13}
\psfrag{14}{14}
\psfrag{15}{15}
\psfrag{16}{16}
\psfrag{17}{17}
\psfrag{18}{18}
\psfrag{19}{19}
\psfrag{20}{20}
\psfrag{21}{21}
\includegraphics[width=.7\linewidth]{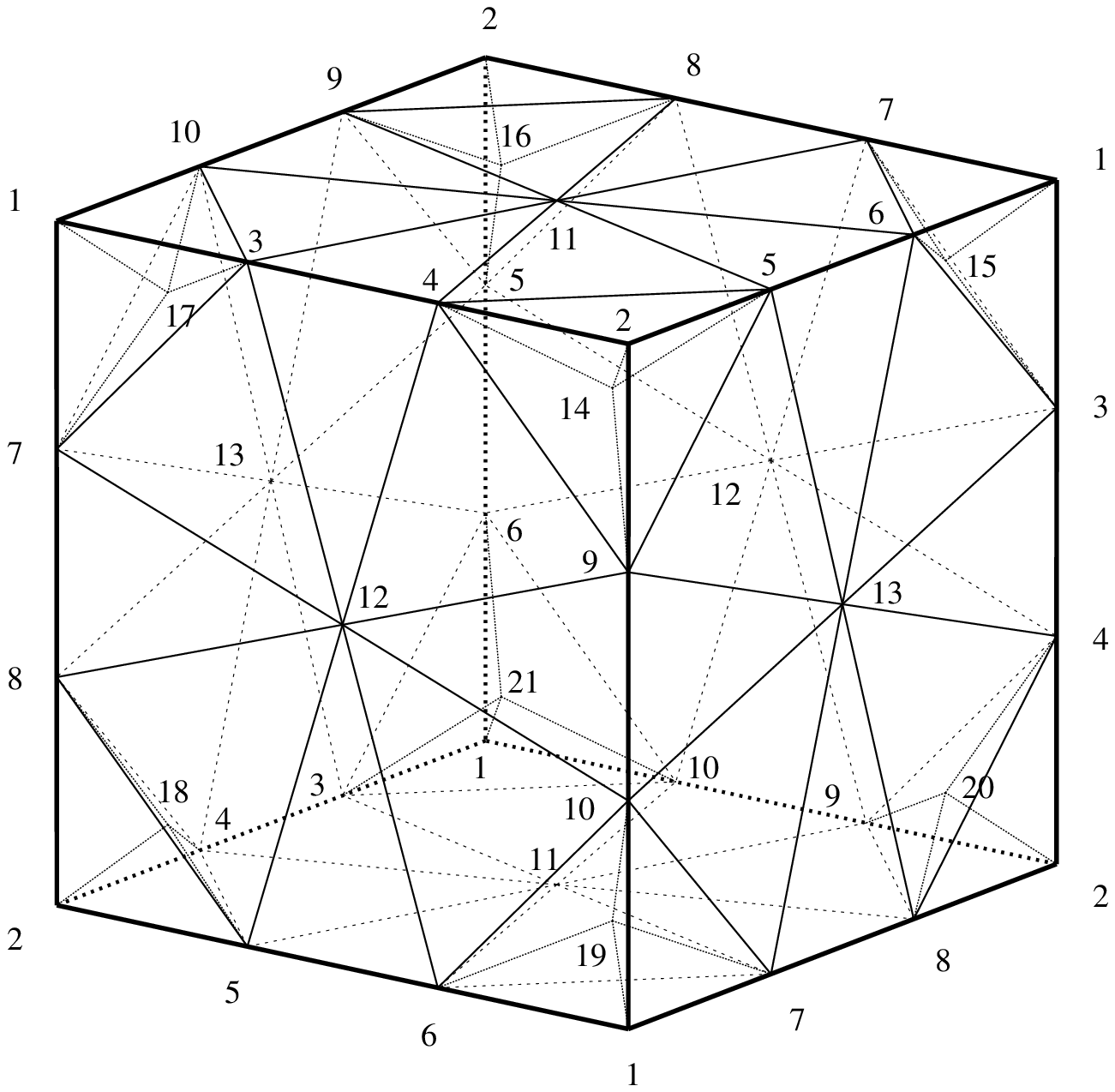}
\end{center}
\caption{The cube or prism space $P(2)$.}
\label{fig:prism_space_2}
\end{figure}

\begin{figure}
\begin{center}
\footnotesize
\psfrag{1}{1}
\psfrag{2}{2}
\psfrag{3}{3}
\psfrag{4}{4}
\psfrag{5}{5}
\psfrag{6}{6}
\psfrag{7}{7}
\psfrag{8}{8}
\psfrag{9}{9}
\psfrag{10}{10}
\psfrag{11}{11}
\psfrag{12}{12}
\psfrag{13}{13}
\psfrag{14}{14}
\psfrag{15}{15}
\psfrag{16}{16}
\psfrag{17}{17}
\psfrag{18}{18}
\psfrag{19}{19}
\includegraphics[width=.825\linewidth]{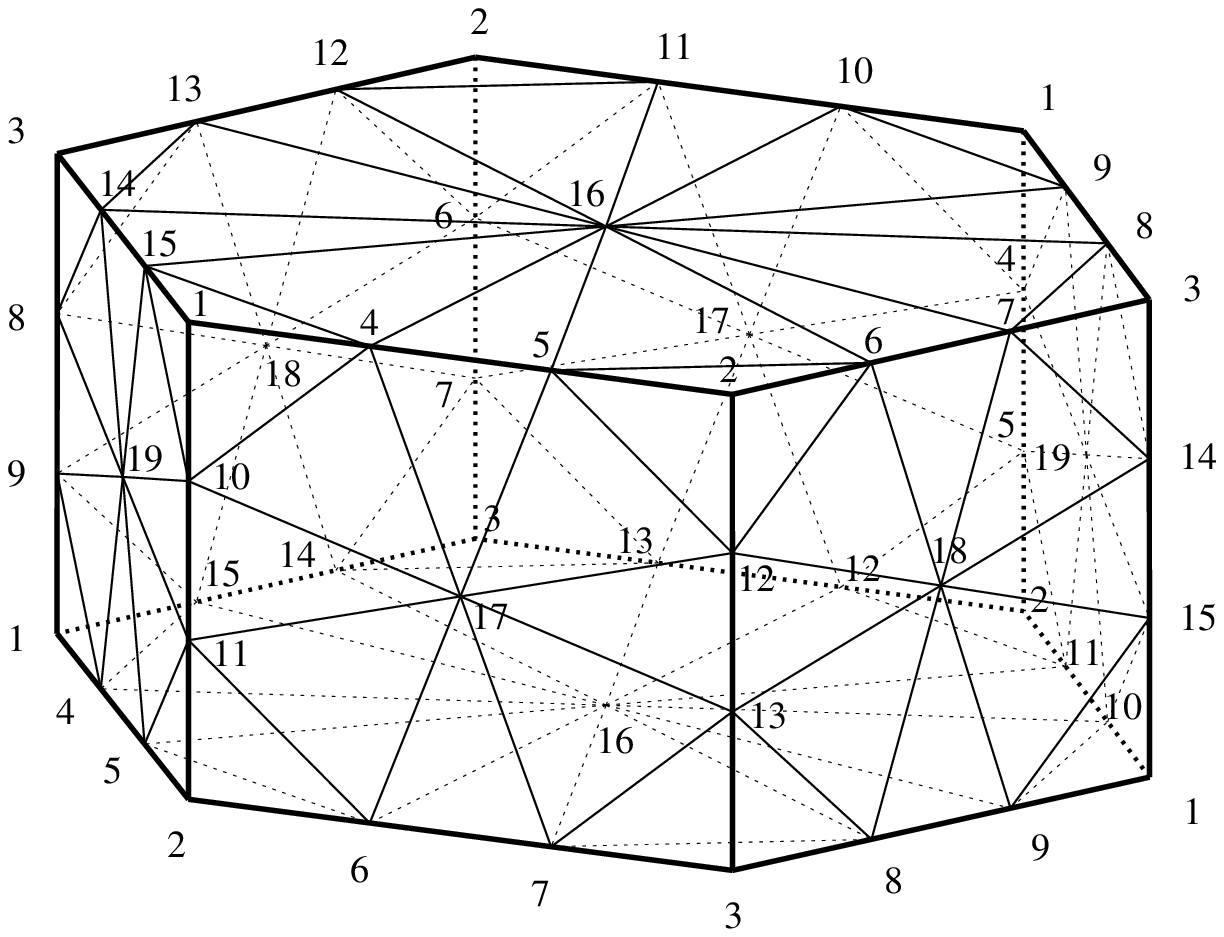}
\end{center}
\caption{The prism space $P(3)$.}
\label{fig:prism_space_3}
\end{figure}

The Seifert manifolds $\{ Oo,0\!\mid\! -1;(2,1),(2,1),(r,1)\}$
are called \emph{prism~spaces} and are denoted by $P(r)$.
Geometrically, $P(r)$ can be obtained from a prism over
a regular $2r$-gon by identifying the top and the bottom side 
under a twist of $\pi/r$. Opposite square side faces are identified
after a twist of $\pi/2$; see Figures~\ref{fig:prism_space_2}
and \ref{fig:prism_space_3}. 

For triangulating the prism spaces, we proceed similar to
\cite{BjoernerLutz2000}.
We first triangulate the boundary of the solid prism by 
introducing a midpoint for every side face
and by subdividing every edge with two additional vertices.
Under the identification, the original $4r$ vertices of the prism
are mapped to $r$ vertices ($1,\ldots,r$ in Figures~\ref{fig:prism_space_2}
and \ref{fig:prism_space_3}). In order to avoid identifications of interior edges,
these vertices are cut off by placing a new vertex just below them.
Then we put an upright bipyramid over a $2r$-gon into the interior of 
the prism (with a vertex of the $2r$-gon of the bipyramid
corresponding to every vertical side face of the prism).
We triangulate the space between the bipyramid and the prism 
(with cut off vertices $1,\ldots,r$) consistently as described in
\cite{BjoernerLutz2000}.
The bipyramid itself is triangulated
by connecting its top with its bottom vertex with an edge
and then slicing the bipyramid around this edge into tetrahedra.

\begin{thm}
There is a series\, $P_{12r+3}(r)$\, of\, ${\mathbb Z}_{\,2r}$-invariant triangulations
of the prism spaces $P(r)$ with\, $f(P_{12r+3}(r))=(12r+3,84r+3,144r,72r)$.
\end{thm}

The prism spaces $P(r)$ have homology groups
$$H_*(P(r))=\left\{
           \begin{array}{l@{\hspace{5mm}}l}
              ({\mathbb Z},{\mathbb Z}_2^{\,2},0,{\mathbb Z}),&\mbox{\rm $r$ even,}\\[1mm]
              ({\mathbb Z},{\mathbb Z}_4,0,{\mathbb Z}),      &\mbox{\rm $r$ odd,}
           \end{array}
           \right.$$
and fundamental group $D^*_{r}$. They arise as quotients $S^3/D^*_{r}$
of the Lie group $S^3$ with respect to left/right 
multiplication of the finite subgroup $D^*_{r}$ of $S^3$.

From the series $P_{12r+3}(r)$ we obtained much smaller triangulations
via bistellar flips; see Table~\ref{tbl:small_prism_spaces}.

\begin{table}
\small\centering
\defaultaddspace=0.3em
\caption{Smallest known triangulations of prism spaces $P(r)$.}\label{tbl:small_prism_spaces}
\begin{tabular}{@{}l@{\hspace{10mm}}l@{}}
\\
\toprule
 \addlinespace
    Manifold         & $f$-Vector  \\ 
\midrule
 \addlinespace
 \addlinespace
 $P(2)=S^3/Q$        & (15,90,150,75) \\
 \addlinespace
 $P(3)$              & (15,97,164,82) \\
 \addlinespace
 $P(4)$              & (15,104,178,89) \\
 \addlinespace
 $P(5)$              & (17,123,212,106) \\
 \addlinespace
 $P(6)$              & (18,135,234,117) \\
 \addlinespace
 $P(7)$              & (19,148,258,129) \\
 \addlinespace
 $P(8)$              & (19,156,274,137) \\
 \addlinespace
 $P(9)$              & (20,169,298,149) \\
 \addlinespace
 $P(10)$             & (21,182,322,161) \\
 \addlinespace
\bottomrule
\end{tabular}
\end{table}

Besides the prism spaces $\{ Oo,0\!\mid\! -1;(2,1),(2,1),(r,1)\}$,
there are further spherical $3$-manifolds corresponding to
the platonic triple $(2,2,\alpha_3)$, namely
the \emph{generalized prism spaces}\, $\{ Oo,0\mid b;(2,1),(2,1),(\alpha_3,\beta_3)\}$.
They are the only spherical $3$-manifolds that allow a fibration $\{ OnI,1\mid b;(\alpha_1,\beta_1)\}$
over the nonorientable orbit surface ${\mathbb R}{\bf P}^{\,2}$.

\bigskip

Associated with the platonic triple $(2,3,5)$ are the spherical dode\-ca\-he\-dron space $S^3/I^*$ 
with Seifert invariants\, $\{ Oo,0\!\mid\! -1;(2,1),(3,1),(5,1)\}$\,
and the \emph{generalized spherical dodecahedral spaces}\, 
$\{ Oo,0\!\mid\! b;(2,1),(3,\beta_2),(5,\beta_3)\}$
(see Table~\ref{tbl:ssm} for the fundamental
groups of these spaces).

The spherical dodecahedral space, i.e., the Poincar\'{e} homology $3$-sphere, 
is the only fibered homology $3$-sphere 
with finite fundamental group \cite{Seifert1933}.
It coincides with the Brieskorn homology sphere $\Sigma (2,3,5)$;
see Section~\ref{sec:hypersurfaces}. Triangulations 
of the Poincar\'{e} homology $3$-sphere,
in particular, a triangulation $\Sigma^{\,3}_{16}$ with $16$ vertices, are given in
\cite{BjoernerLutz2000}; see also \cite{BjoernerLutz2003}.
Triangulations of the generalized spherical dodecahedral spaces 
are constructed in \cite{BrehmLutz2002pre}.

\bigskip

The spherical spaces that are obtained for the platonic triple $(2,3,3)$
are the octahedral space $S^3/T^*$, with the binary tetrahedral group $T^*$  
of order $24$ as fundamental group, and the generalized octahedral spaces.

\begin{figure}
\begin{center}
\footnotesize
\psfrag{1}{1}
\psfrag{2}{2}
\psfrag{3}{3}
\psfrag{4}{4}
\psfrag{5}{5}
\psfrag{6}{6}
\psfrag{7}{7}
\psfrag{8}{8}
\psfrag{9}{9}
\psfrag{10}{10}
\psfrag{11}{11}
\psfrag{12}{12}
\psfrag{13}{13}
\psfrag{14}{14}
\psfrag{15}{15}
\psfrag{16}{16}
\psfrag{17}{17}
\psfrag{18}{18}
\psfrag{19}{19}
\includegraphics[width=.865\linewidth]{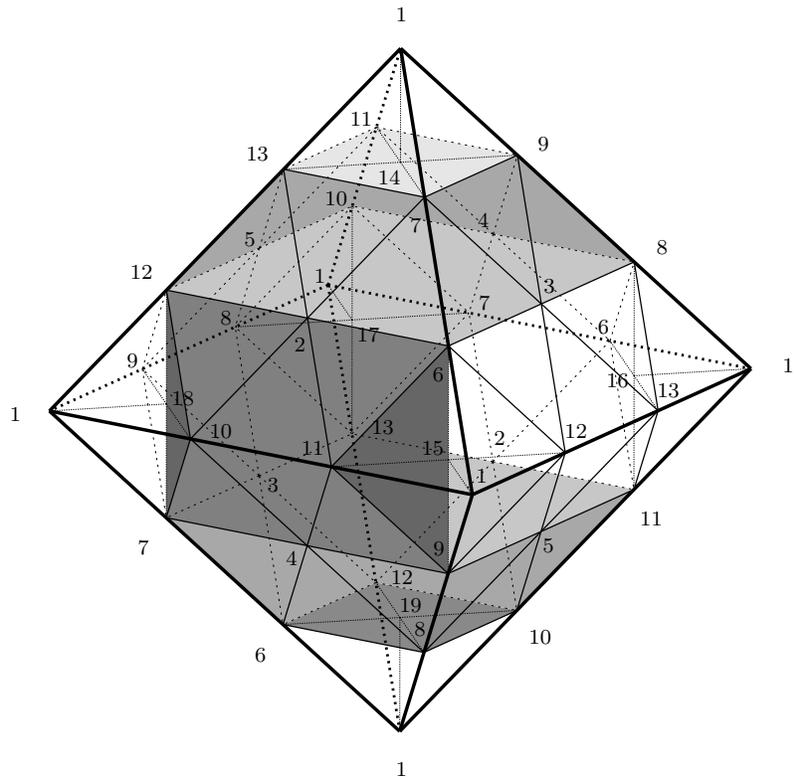}
\end{center}
\caption{Octahedral space.}
\label{fig:octa_space}
\end{figure}

The octahedral space can be triangulated as follows. We begin with a solid
octahedron and identify opposite triangles by a twist of $\pi/3$.
We subdivide every edge with two vertices and every triangle 
by inserting a center vertex; see Figure~\ref{fig:octa_space}. 
Under the identification, all original vertices are mapped to the vertex $1$
and are cut off by placing new vertices below them.
After the cutting, the boundary (with identifications) of the 
truncated octahedron is split into two polar and four equatorial
regions as depicted in Figure~\ref{fig:octa_space}.
We place another octahedron in the center of the truncated octahedron,
with one equatorial respectively polar vertex of the inner octahedron corresponding
to every equatorial respectively polar region of the outer truncated octahedron.
The space between the outer truncated and the inner octahedron is triangulated
consistently. The inner octahedron is triangulated by introducing an
edge from top to bottom.

The remaining platonic triple $(2,3,4)$ yields the truncated cube
space $S^3/O^*$ and the generalized truncated cube spaces; 
see Table~\ref{tbl:ssm}.

\begin{figure}
\begin{center}
\footnotesize
\psfrag{1}{1}
\psfrag{2}{2}
\psfrag{3}{3}
\psfrag{4}{4}
\psfrag{5}{5}
\psfrag{6}{6}
\psfrag{7}{7}
\psfrag{8}{8}
\psfrag{9}{9}
\psfrag{10}{10}
\psfrag{11}{11}
\psfrag{12}{12}
\psfrag{13}{13}
\psfrag{14}{14}
\psfrag{15}{15}
\psfrag{16}{16}
\includegraphics[width=.79\linewidth]{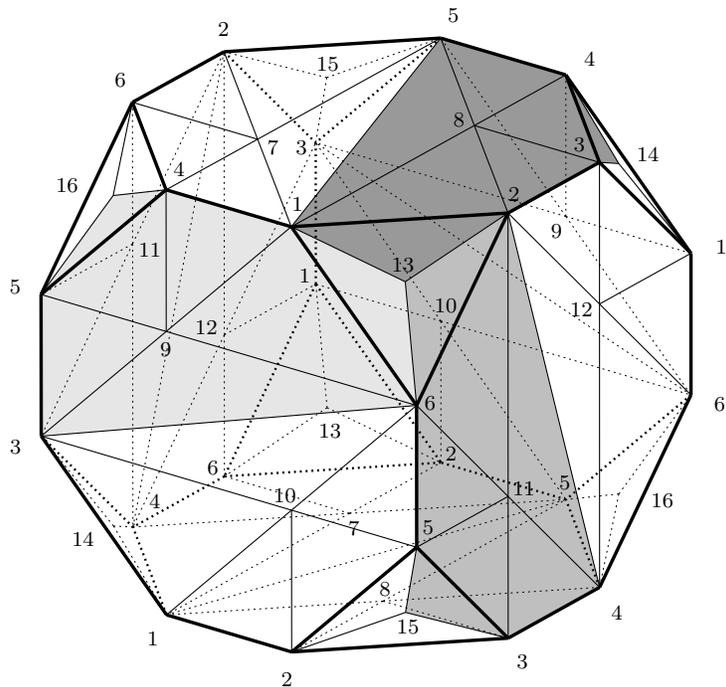}
\end{center}
\caption{Truncated cube space.}
\label{fig:truncated_cube}
\end{figure}

For a triangulation of $S^3/O^*$ we triangulate first 
the boundary of a truncated cube with identified opposite
octagons and triangles as in Figure~\ref{fig:truncated_cube}.
Then we split the (identified) boundary into $12$ congruent 
pieces (three of them are shaded in Figure~\ref{fig:truncated_cube}).
To every of the $12$ pieces we associate a vertex of an icosahedron
that we place into the interior of the truncated cube and then triangulate
the space consistently.

Once again, by applying bistellar flips to the above triangulations
of the spherical octahedron and the truncated cube space, we obtain
small triangulations of these spaces. For the $f$-vectors of these 
triangulations see Table~\ref{tbl:small_platonic_spaces}.
The full triangulations can be found online at \cite{PAGE}.

\begin{thm}
The spherical octahedron and the truncated cube space can be
triangulated with $15$ respectively $16$ vertices.
\end{thm}

\begin{conj}
The triangulations with $15$ respectively $16$ vertices
of the spherical octahedron and the truncated cube space
are vertex-minimal.
\end{conj}

\begin{table}
\small\centering
\defaultaddspace=0.3em
\caption{Smallest known triangulations of the spherical octahedral,
  truncated cube and dodecahedral space.}\label{tbl:small_platonic_spaces}
\begin{tabular}{@{}l@{\hspace{10mm}}l@{\hspace{10mm}}l@{}}
\\
\toprule
 \addlinespace
    Manifold         & Homology & $f$-Vector \\ 
\midrule
 \addlinespace
 \addlinespace
 $S^3/T^*$           & $({\mathbb Z},{\mathbb Z}_3,0,{\mathbb Z})$ & (15,102,174,87) \\
 \addlinespace
 $S^3/O^*$           & $({\mathbb Z},{\mathbb Z}_2,0,{\mathbb Z})$ & (16,109,186,93) \\
 \addlinespace
 $\Sigma^{\,3}=S^3/I^*$ & $({\mathbb Z},0,0,{\mathbb Z})$& (16,106,180,90) \\
 \addlinespace
\bottomrule
\end{tabular}
\end{table}

Now that we have seen various examples of spherical $3$-manifolds,
we return to Thurston's Geometrization Conjecture and its implications
for $3$-manifolds with finite fundamental groups.
\begin{conj} {\rm (Thurston's Elliptization Conjecture \cite{Thurston1982})}
Every (closed) $3$-manifold with finite fundamental group can be modeled
on the spherical geometry $(S^3,SO(4))$.
\end{conj}
If a $3$-manifold has finite fundamental group, then its universal cover is a 
homotopy sphere. Thus, in particular, Thurston's Elliptization
Conjecture comprises the Poincar\'e Conjecture. 
(A proof of the Elliptization Conjecture 
has recently been announced 
by Perelman \cite{Perelman2003bpre}.)

Thurston's conjecture also states that every free action 
of a finite group on~$S^3$ is conjugate to an orthogonal action,
and, at least for this part of the conjecture, partial results 
have been known for a while. Milnor \cite{Milnor1957} determined a 
list of finite groups which might act freely on $S^3$. 
This list was reduced by Lee \cite{Lee1973},
and (based on the Smale Conjecture that $\text{Diff}(S^3)=O(4)$, 
which later was proved by Hatcher \cite{Hatcher1983}) 
Thomas \cite{Thomas1978} showed that, in fact, only the subgroups 
of $SO(4)$ classified by Hopf and Threlfall and Seifert 
need to be considered. For a number of these groups
it has been shown that if they act fixed point freely on $S^3$, 
then this action is conjugate to an orthogonal action.
See Myers \cite{Myers1981}, Scott
\cite{Scott1983}, and Maher and Rubinstein \cite{MaherRubinstein2003}
for a discussion and for further references.


\section{Flat Spaces}

Discrete groups of isometries of the Euclidean space $E^3$
were first studied in the context of crystallography 
in the late 19th century (see the references in \cite{Wolf1967}). 
In 1911, Bieberbach (\cite{Bieberbach1911}, \cite{Bieberbach1912}) 
presented a comprehensive structure theory for crystallographic groups. 
Based on the results by Bieberbach, Nowacki \cite{Nowacki1934} 
gave a homeomorphism classification and Hantzsche and Wendt 
\cite{HantzscheWendt1935} worked out an affine classification 
of Euclidean space forms. An isometric classification of flat $3$-manifolds
was later obtained by Wolf \cite[Ch.~3]{Wolf1967}.

\begin{thm} \rm{(Nowacki \cite{Nowacki1934}, Hantzsche and Wendt 
\cite{HantzscheWendt1935}, and Wolf \cite{Wolf1967})}
There are exactly $6$ orientable flat $3$-manifolds, 
${\cal G}_1$, \ldots, ${\cal G}_6$, and $4$ non-orientable ones,
${\cal B}_1$, \ldots, ${\cal B}_4$.
\end{thm}

Nine of the ten flat spaces can be obtained as torus or Klein bottle
bundles over $S^1$. The direct products $T^2\!\times S^1=T^3$
and $K\!\times S^1$ yield the manifolds ${\cal G}_1$ and ${\cal B}_1$,
respectively. The other seven bundles are twisted products.
All these spaces can easily be triangulated. 

\bigskip

\begin{figure}
\begin{center}
\footnotesize
\psfrag{1}{1}
\psfrag{2}{2}
\psfrag{3}{3}
\psfrag{4}{4}
\psfrag{5}{5}
\psfrag{6}{6}
\psfrag{7}{7}
  \includegraphics[width=.625\linewidth]{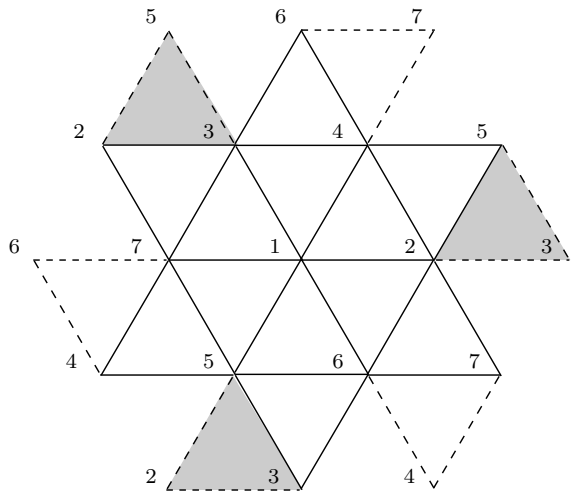}
\end{center}
\caption{M\"obius torus.}
\label{fig:torus7}
\end{figure}

\emph{Triangulations of orientable flat spaces}.
The unique minimal triangulation of the $2$-torus $T^2$ is M\"obius' torus \cite{Moebius1886}
with $7$ vertices and is depicted in Figure~\ref{fig:torus7}.

The direct product triangulation of the M\"obius torus with
$S^1$ triangulated as an empty triangle has $7\cdot 3=21$ vertices.
A $15$-vertex triangulation of $T^3$ was first found by 
K\"uhnel and Lassmann \cite{KuehnelLassmann1984-3torus},
and it is conjectured 
that $15$ vertices are minimal for a triangulation of the $3$-dimensional torus.

M\"obius' torus has the affine group $AGL(1,7)$ as vertex-transitive combinatorial 
automorphism group, and every group element of $AGL(1,7)$ can be used to construct
a corresponding twisted product triangulation of the $2$-torus over 
$S^1$ by applying this group element to one boundary component 
of the product triangulation of\, $T^2\!\times I$\, before gluing
together the two boundary $2$-tori. 
Gluing under the permutations $(2,4,3,7,5,6)$, $(2,3,5)(4,7,6)$, and $(2,7)(3,6)(4,5)$
of the vertices of the M\"obius torus in Figure~\ref{fig:torus7}, 
which induce rotations by $2\pi/6$, $2\pi/3$, and $2\pi/2$, 
yields the flat spaces ${\cal G}_5$, ${\cal G}_3$, and ${\cal G}_2$, respectively.
All other twisted products of the M\"obius torus over $S^1$
are homeomorphic to one of these spaces (or to $T^3={\cal G}_1$ under
the trivial identification).

\begin{figure}
\begin{center}
\footnotesize
\psfrag{1}{1}
\psfrag{2}{2}
\psfrag{3}{3}
\psfrag{4}{4}
\psfrag{5}{5}
\psfrag{6}{6}
\psfrag{7}{7}
\psfrag{8}{8}
\psfrag{9}{9}
\psfrag{10}{10}
\includegraphics[width=.47\linewidth]{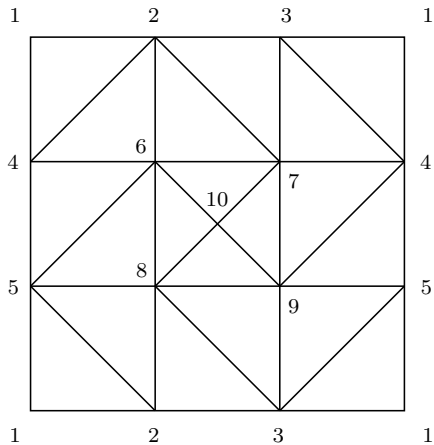}
\end{center}
\caption{Torus with rotation by $\pi/2$.}
\label{fig:torus18}
\end{figure}

There is one further orientable flat $T^2$-bundle over $S^1$,
the manifold ${\cal G}_4$, which can be triangulated
by starting with the $10$-vertex triangulation of the $2$-torus
in Figure~\ref{fig:torus18} and then identifying the two boundary components of
the corresponding twisted product triangulation over $I$
under a twist of $2\pi/4$. 
In fact, there is no smaller twisted product triangulation of ${\cal G}_4$. 
There are $7$ combinatorially distinct triangulations of $T^2$
with $8$ vertices \cite{DattaNilakantan2002} and $112$ triangulations
with $9$ vertices \cite{Lutz2003dpre},
but none of these admits a (combinatorial) rotation by $\pi/2$.
\begin{thm}
The orientable flat spaces ${\cal G}_2$, ${\cal G}_3$, ${\cal G}_4$, and ${\cal G}_5$
have minimal twisted product triangulations with $21$, $21$, $30$, and $21$
vertices, respectively. 
\end{thm}

The remaining orientable flat $3$-manifold ${\cal G}_6$ is not
a torus bundle over~$S^1$, but can be constructed by gluing together 
two solid cubes along their boundary squares in a very symmetric way,
as described in \cite[p.~126]{Thurston1997}. 
\begin{figure}
\begin{center}
\footnotesize
\psfrag{A}{A}
\psfrag{B}{B}
\psfrag{C}{C}
\psfrag{D}{D}
\psfrag{E}{E}
\psfrag{F}{F}
\psfrag{G}{G}
\psfrag{H}{H}
\psfrag{I}{I}
\psfrag{J}{J}
\psfrag{K}{K}
\psfrag{L}{L}
\psfrag{M}{M}
\psfrag{N}{N}
\psfrag{O}{O}
\psfrag{P}{P}
\psfrag{Q}{Q}
\psfrag{R}{R}
\psfrag{S}{S}
\psfrag{T}{T}
\psfrag{U}{U}
\psfrag{V}{V}
\psfrag{W}{W}
\psfrag{X}{X}
\includegraphics[width=\linewidth]{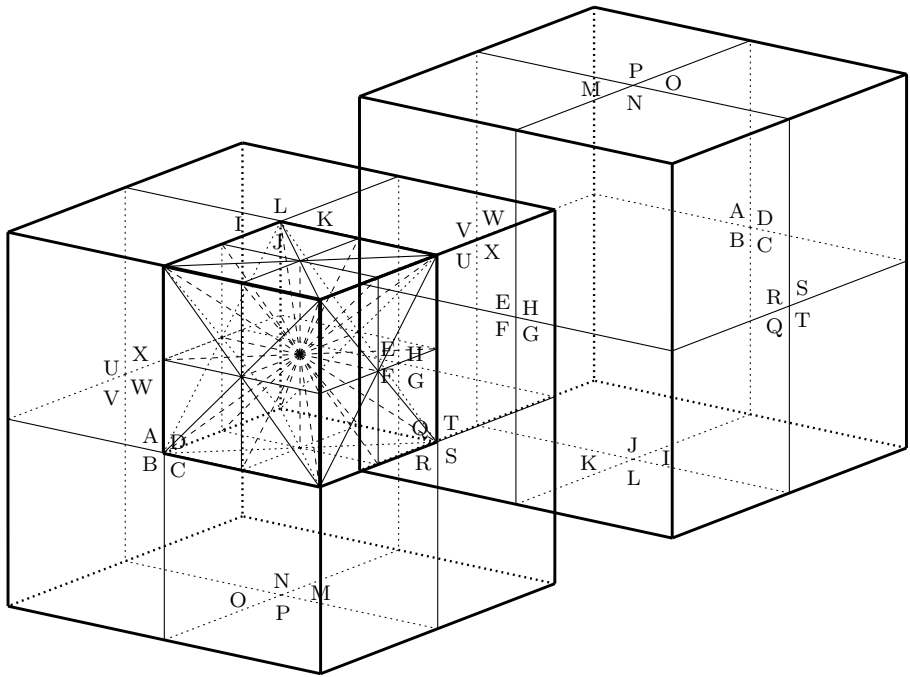}
\end{center}
\caption{The flat space ${\cal G}_6$.}
\label{fig:g6}
\end{figure}
In Figure~\ref{fig:g6}, the two cubes are shown with the respective
identifications. For a triangulation of this space, 
the two cubes are each split into eight cubes, 
which then are triangulated by their barycentric subdivision. 
The resulting triangulation of ${\cal G}_6$ 
has $f$-vector\, $f=(128,896,1536,768)$.  

Since flat $3$-manifolds are Seifert manifolds, 
we can alternatively obtain triangulations of them 
as described in \cite{BrehmLutz2002pre}.
The Euler number $e(M)$ of an Euclidean $3$-manifold $M$ and
the Euler characteristic of the respective base orbifold are both zero.
It is then an easy exercise to compute the sets of Seifert invariants,
which are possible for the flat spaces.

By applying bistellar flips to the above triangulations
we get small triangulations.
\begin{thm}
The orientable flat $3$-manifolds ${\cal G}_1$, \ldots, ${\cal G}_6$
can be triangulated with $15$, $17$, $17$, $16$, $16$, and $17$ vertices, respectively.
\end{thm}
The homology groups, the Seifert fibrations, and the smallest known 
$f$-vectors of the spaces ${\cal G}_1$, \ldots, ${\cal G}_6$ and of
the spaces ${\cal B}_1$, \ldots, ${\cal B}_4$ are listed
in Table~\ref{tbl:small_flat_spaces}. Note that the ten flat spaces
have pairwise distinct homology. The spaces ${\cal G}_2$, ${\cal B}_1$,
${\cal B}_2$, and ${\cal B}_4$ have two different Seifert fibrations,
respectively, while the other spaces each have a unique Seifert fibration.

\bigskip

\begin{figure}
\begin{center}
\footnotesize
\psfrag{1}{1}
\psfrag{2}{2}
\psfrag{3}{3}
\psfrag{4}{4}
\psfrag{5}{5}
\psfrag{6}{6}
\psfrag{7}{7}
\psfrag{8}{8}
\includegraphics[width=.6\linewidth]{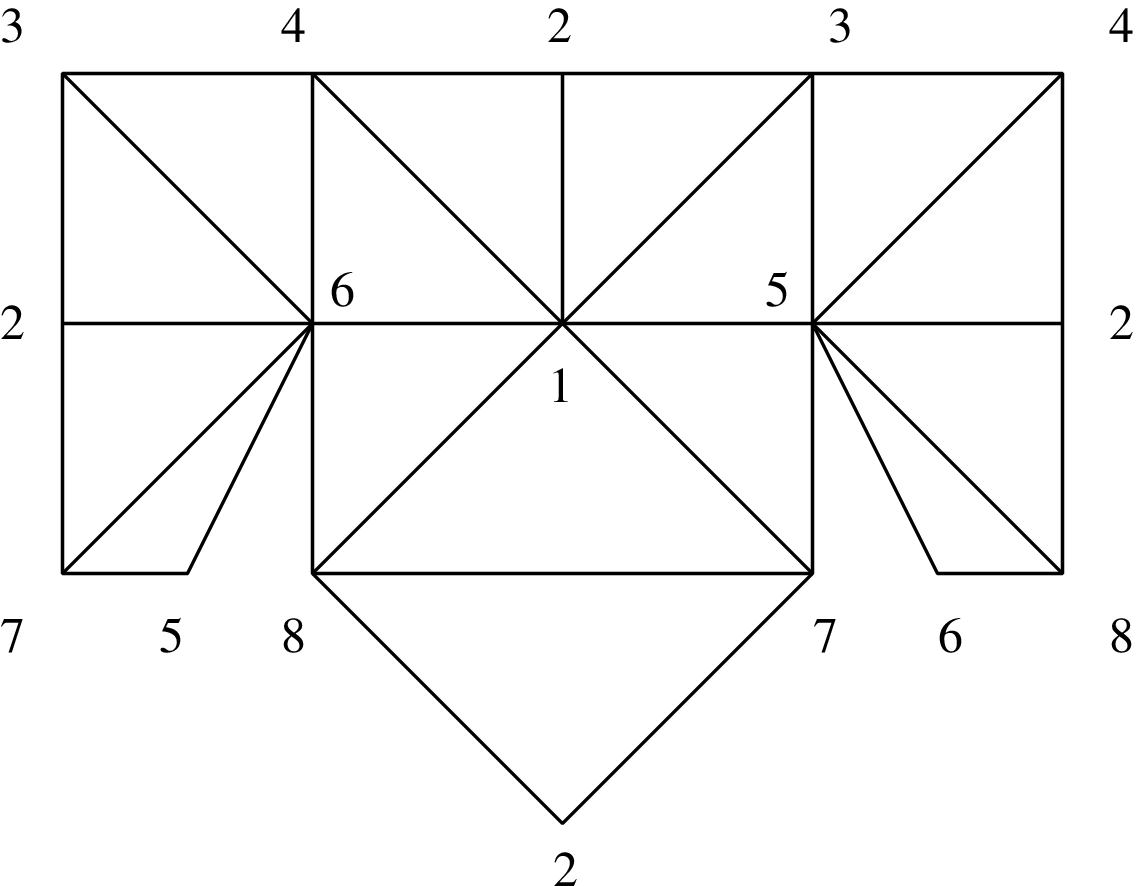}
\end{center}
\caption{$8$-vertex Klein bottle\, $\mbox{}^2\hspace{.3pt}8_{\,20}$.}
\label{fig:Klein8}
\end{figure}

\emph{Triangulations of non-orientable flat spaces}.
All four non-orientable flat $3$-manifolds ${\cal B}_1$, \ldots,
${\cal B}_4$ are Klein bottle bundles over $S^1$.
(The space ${\cal B}_2$ can also be obtained as a non-orientable torus
bundle over $S^1$.) The Klein bottle $K$ can be minimally be triangulated
with $8$ vertices. There are althogether six $8$-vertex triangulations
of $K$ (\cite{Cervone1994}, \cite{DattaNilakantan2002}), one
of which is the example $\mbox{}^2\hspace{.3pt}8_{\,20}$ of \cite{Lutz2003dpre},
shown in Figure~\ref{fig:Klein8}. It can be read off directly from the figure
that the example is a triangulation of the Klein bottle, since the
boundary consists of $2+2$ arcs 2--7--5--6--8--2, 2--7--5--6--8--2,
2--4--3--2, and 2--4--3--2, which are just the arcs of the 
fundamental polygon (see \cite[Ch.~1]{Stillwell1980}) of the Klein
bottle.

The automorphism group of $\mbox{}^2\hspace{.3pt}8_{\,20}$
is isomorphic to the dihedral group $D_4$ and is rich enough
for constructing (twisted) product triangulations 
of all four non-orientable flat $3$-manifolds
from the triangulation $\mbox{}^2\hspace{.3pt}8_{\,20}$.
In Table~\ref{tbl:b1_b4} we list the group elements
and the corresponding flat spaces.
\begin{thm}
The non-orientable flat spaces ${\cal B}_1$, ${\cal B}_2$, ${\cal B}_3$, 
and ${\cal B}_4$ have minimal (twisted) product triangulations with $24$ vertices.
\end{thm}
Taking these triangulations as starting triangulations,
we obtained smaller triangulations via bistellar flips.
\begin{thm}
The non-orientable flat $3$-manifolds ${\cal B}_1$, \ldots, ${\cal B}_4$
can be triangulated with $17$, $16$, $17$, and $17$ vertices, respectively.
\end{thm}

\begin{table}
\small\centering
\defaultaddspace=0.3em
\caption{Twisted product spaces of the Klein bottle\, $\mbox{}^2\hspace{.3pt}8_{\,20}$.}\label{tbl:b1_b4}
\begin{tabular}{@{}l@{\hspace{10mm}}l@{}}
\\
\toprule
 \addlinespace
    Manifold   & Gluing Transformations \\ 
\midrule
 \addlinespace
 \addlinespace
 ${\cal B}_1$  & $()$, \\
               & $(1,2)(5,6)$ \\
 \addlinespace
 ${\cal B}_2$  & $(1,5,2,6)(3,8)(4,7)$, \\
               & $(1,6,2,5)(3,8)(4,7)$ \\
 \addlinespace
 ${\cal B}_3$  & $(1,2)(3,4)(7,8)$, \\
               & $(3,4)(5,6)(7,8)$ \\
 \addlinespace
 ${\cal B}_4$  & $(1,5)(2,6)(3,7)(4,8)$,  \\
               & $(1,6)(2,5)(3,7)(4,8)$ \\
 \addlinespace
\bottomrule
\end{tabular}
\end{table}

In this section we have chosen to construct triangulations of the
flat spaces via products and twisted products (except for ${\cal G}_6$).
This approach has the advantage that the resulting triangulations
are rather small and simple to produce. Identifying the 
faces of a cube, of a hexagonal prism, 
or of a rhombic dodecahedron, however, 
gives a more geometric way for obtaining flat spaces.
See Cipra \cite{Cipra2002a} for illustrations of the identifications.

\begin{landscape}

\begin{table}
\small\centering
\defaultaddspace=0.3em
\caption{Smallest known triangulations of flat spaces.}\label{tbl:small_flat_spaces}
\begin{tabular*}{\linewidth}{@{\extracolsep{\fill}}llll@{}}
\\
\toprule
 \addlinespace
    Manifold         & Seifert Fibrations & Homology & $f$-Vector \\
\midrule
 \addlinespace
 \addlinespace
 ${\cal G}_1=T^3$    & $\{ Oo,1\mid 0\}$ & $({\mathbb Z},{\mathbb Z}^3,{\mathbb Z}^3,{\mathbb Z})$ & (15,105,180,90) \\
 \addlinespace
 ${\cal G}_2$ & $\{ Oo,0\mid -2;(2,1),(2,1),(2,1),(2,1)\}$, & $({\mathbb Z},{\mathbb Z}\oplus{\mathbb Z}_2^{\,2},{\mathbb Z},{\mathbb Z})$ & (17,118,202,101) \\
                                                                                            & $\{ On,2\mid 0\}$ & & \\
 \addlinespace
 ${\cal G}_3$ & $\{ Oo,0\mid -1;(3,1),(3,1),(3,1)\}$ & $({\mathbb Z},{\mathbb Z}\oplus{\mathbb Z}_3,{\mathbb Z},{\mathbb Z})$ & (17,117,200,100) \\
 \addlinespace
 ${\cal G}_4$ & $\{ Oo,0\mid -1;(2,1),(4,1),(4,1)\}$ & $({\mathbb Z},{\mathbb Z}\oplus{\mathbb Z}_2,{\mathbb Z},{\mathbb Z})$ & (16,115,198,99) \\
 \addlinespace
 ${\cal G}_5$ & $\{ Oo,0\mid -1;(2,1),(3,1),(6,1)\}$ & $({\mathbb Z},{\mathbb Z},{\mathbb Z},{\mathbb Z})$ & (16,112,192,96) \\
 \addlinespace
 ${\cal G}_6$        & $\{ On,1\mid -1;(2,1),(2,1)\}$ & $({\mathbb Z},{\mathbb Z}_4^{\,2},0,{\mathbb Z})$ & (17,124,214,107) \\
 \addlinespace
\midrule
 \addlinespace
 \addlinespace
 ${\cal B}_1=K\!\times S^1$ & $\{ NnI,2\mid 0\}$, $\{ No,1\mid 0\}$ & $({\mathbb Z},{\mathbb Z}^2\oplus{\mathbb Z}_2,{\mathbb Z}\oplus{\mathbb Z}_2,0)$ & (17,118,202,101) \\
 \addlinespace
 ${\cal B}_2$ & $\{ NnI,2\mid 1\}$, $\{ No,1\mid 1\}$ & $({\mathbb Z},{\mathbb Z}^2,{\mathbb Z}\oplus{\mathbb Z}_2,0)$ & (16,110,188,94) \\
 \addlinespace
 ${\cal B}_3$ & $\{ NnI\!I,2\mid 0\}$ & $({\mathbb Z},{\mathbb Z}\oplus{\mathbb Z}_2^{\,2},{\mathbb Z}_2,0)$ & (17,119,204,102) \\
 \addlinespace
 ${\cal B}_4$ & $\{ NnI\!I,2\mid 1\}$, $\{ NnI,1\mid 0;(2,1),(2,1)\}$ & $({\mathbb Z},{\mathbb Z}\oplus{\mathbb Z}_4,{\mathbb Z}_2,0)$ & (17,117,200,100) \\
 \addlinespace
\bottomrule
\end{tabular*}
\end{table}

\end{landscape}

\section{${\rm Nil}$-Spaces}
\label{sec:nil_spaces}

The Heisenberg group ${\rm Nil}$ of upper triangular matrices $A\in GL(3,\mathbb{R})$\, of the form
\addtolength{\arraycolsep}{.85mm}
$$A=\left(\begin{array}{rrr}1&u&w\\0&1&v\\0&0&1\end{array}\right)\addtolength{\arraycolsep}{-.85mm}$$
is a nilpotent Lie group. The ${\rm Nil}$-spaces modeled on the geometry ${\rm Nil}$
are orientable and can be Seifert fibered 
with Euler number $e\neq 0$ and
base orbifold of Euler characteristic $\chi =0$.
The possible sets of Seifert invariants are 
listed in Table~\ref{tbl:nil_spaces}.


\begin{table}
\small\centering
\defaultaddspace=0.3em
\caption{Seifert fibrations of ${\rm Nil}$-spaces.}\label{tbl:nil_spaces}
\begin{tabular*}{\linewidth}{@{\extracolsep{\fill}}ll@{}}
\\
\toprule
 \addlinespace
Seifert Fibration & Homology  \\ 
\midrule
 \addlinespace
 \addlinespace
$\{ Oo,0\mid b;(2,1),(3,1),(6,1)\}$,\quad $b\neq -1$ & $({\mathbb Z},{\mathbb Z}_{\,36|b+1|},0,{\mathbb Z})$ \\
 \addlinespace
$\{ Oo,0\mid b;(2,1),(3,1),(6,5)\}$                  & $({\mathbb Z},{\mathbb Z}_{\,24+36|b+1|},0,{\mathbb Z})$,\quad $b\geq -1$, \\
                                                     & $({\mathbb Z},{\mathbb Z}_{\,12+36|b+2|},0,{\mathbb Z})$,\quad $b\leq -2$ \\
 \addlinespace
$\{ Oo,0\mid b;(2,1),(4,1),(4,1)\}$,\quad $b\neq -1$ & $({\mathbb Z},{\mathbb Z}_2\oplus {\mathbb Z}_{\,16|b+1|},0,{\mathbb Z})$ \\
 \addlinespace
$\{ Oo,0\mid b;(2,1),(4,1),(4,3)\}$                  & $({\mathbb Z},{\mathbb Z}_2\oplus {\mathbb Z}_{\,16|b+\frac{3}{2}|},0,{\mathbb Z})$ \\
 \addlinespace
$\{ Oo,0\mid b;(3,1),(3,1),(3,1)\}$,\quad $b\neq -1$ & $({\mathbb Z},{\mathbb Z}_3\oplus {\mathbb Z}_{\,9|b+1|},0,{\mathbb Z})$ \\
 \addlinespace
$\{ Oo,0\mid b;(3,1),(3,1),(3,2)\}$                  & $({\mathbb Z},{\mathbb Z}_3\oplus {\mathbb Z}_{\,3+9|b+1|},0,{\mathbb Z})$,\quad $b\geq -1$, \\
                                                     & $({\mathbb Z},{\mathbb Z}_3\oplus {\mathbb Z}_{\,6+9|b+2|},0,{\mathbb Z})$,\quad $b\leq -2$ \\
 \addlinespace
$\{ Oo,0\mid b;(2,1),(2,1),(2,1),(2,1)\}$,\quad $b\neq -2$ & $({\mathbb Z},{\mathbb Z}_2^{\,2}\oplus {\mathbb Z}_{\,4|b+2|},0,{\mathbb Z})$ \\
 \addlinespace
$\{ Oo,\hspace{1.1pt}1\mid b\}$,\quad $b\neq 0$      & $({\mathbb Z},{\mathbb Z}^2\oplus {\mathbb Z}_b,{\mathbb Z}^2,{\mathbb Z})$ \\
 \addlinespace
$\{ On,1\mid b;(2,1),(2,1)\}$,\quad $b\neq -1$       & $({\mathbb Z},{\mathbb Z}_4^{\,2},0,{\mathbb Z})$\\
 \addlinespace
$\{ On,2\mid b\}$,\quad $b\neq 0$                    & $({\mathbb Z},{\mathbb Z}\oplus {\mathbb Z}_2^{\,2},{\mathbb Z},{\mathbb Z})$\quad if\, $b$\, is even, \\
                                                     & $({\mathbb Z},{\mathbb Z}\oplus {\mathbb Z}_4,{\mathbb Z},{\mathbb Z})$\quad if\, $b$\, is odd \\
 \addlinespace
\bottomrule
\end{tabular*}
\end{table}


The ${\rm Nil}$-spaces\, $M_b=\{ Oo,1\mid b\}$\, and\, $M'_b=\{ On,2\mid b\}$
with $b\neq 0$ are examples of ${\mathbb T}^2$-bundles over $S^1$. 
The manifolds $M_b$ also arise as quotients of ${\rm Nil}$ with respect to 
the discrete subgroup of matrices $N_b$ of ${\rm Nil}$ for which $u,v,w\in b{\mathbb Z}$;
see \cite[p.~124]{Orlik1972} and \cite{OrlikRaymond1969}.

In fact, ${\mathbb T}^2$-bundles over $S^1$ are either flat or belong
to the geometries ${\rm Nil}$ or ${\rm Sol}$, depending only on the gluing map; see \cite{Scott1983}
for details.


\mathversion{bold}
\section{($H^2\times{\mathbb R}$)-Spaces}
\mathversion{normal}

Products of hyperbolic surfaces with $S^1$ are natural examples of ($H^2\times{\mathbb R}$)-spaces.
Also the mapping tori of isometries of hyperbolic surfaces are ($H^2\times{\mathbb R}$)-spaces.

For triangulating the product spaces we can start with any
triangulation of an orientable surface $M^2_{(+,g)}$ of
genus $g\geq 2$ or a non-orient\-able surface $M^2_{(-,g)}$ of genus $g\geq 3$
and then form the product triangulation. Vertex-minimal triangulations
of surfaces were mainly obtained by Ringel \cite{Ringel1955} and 
Jungerman and Ringel \cite{JungermanRingel1980}; see also the book by
Ringel \cite{Ringel1974}. The $f$-vectors of small triangulations of 
some products can be found in Table~\ref{tbl:small_HxR_spaces}.
The product spaces
have homology\,
$H_*(M^2_{(+,g)}\!\times\!S^1)=({\mathbb Z},{\mathbb Z}^{\,2g+1},{\mathbb Z}^{\,2g+1},{\mathbb Z})$\,
and\,
$H_*(M^2_{(-,g)}\!\times\!S^1)=({\mathbb Z},{\mathbb Z}^g\oplus {\mathbb Z}_2,{\mathbb Z}^{\,g-1}\oplus {\mathbb Z}_2,0)$.

\begin{table}
\small\centering
\defaultaddspace=0.3em
\caption{Smallest known triangulations of some ($H^2\times{\mathbb R}$)-spaces.}\label{tbl:small_HxR_spaces}
\begin{tabular}{@{}l@{\hspace{10mm}}l@{}}
\\
\toprule
 \addlinespace
    Manifold         & $f$-Vector \\ 
\midrule
 \addlinespace
 \addlinespace
 $M^2_{(+,2)}\!\times\!S^1$    & (20,171,302,151) \\
 \addlinespace
 $M^2_{(+,3)}\!\times\!S^1$    & (22,211,378,189) \\
 \addlinespace
 $M^2_{(+,4)}\!\times\!S^1$    & (26,269,486,243) \\
 \addlinespace
 $M^2_{(+,5)}\!\times\!S^1$    & (28,312,568,284) \\
 \addlinespace
\midrule
 \addlinespace
 \addlinespace
 $M^2_{(-,3)}\!\times\!S^1$    & (19,146,254,127) \\
 \addlinespace
 $M^2_{(-,4)}\!\times\!S^1$    & (20,168,296,148) \\
 \addlinespace
 $M^2_{(-,5)}\!\times\!S^1$    & (21,192,342,171) \\
 \addlinespace
\bottomrule
\end{tabular}
\end{table}

An abundance of twisted product triangulations of ($H^2\times{\mathbb R}$)-spaces
can be constructed with the help of combinatorial isomorphisms of
triangulated hyperbolic surfaces.

\mathversion{bold}
\section{$\widetilde{SL}(2,{\mathbb R})$-Spaces}
\mathversion{normal}

The last of the six Seifert geometries is $\widetilde{SL}(2,{\mathbb R})$
and is discussed in \cite{Scott1983}, \cite{Thurston1997}, and \cite{Thurston2002}. 
Particularly interesting examples of $\widetilde{SL}(2,{\mathbb R})$-spaces are 
the Seifert homology spheres, which, with the exception of the
Poincar\'e sphere with finite fundamental group, all belong to this geometry.
See Section~\ref{sec:homology_spheres} for further comments.

\section{${\rm Sol}$-Spaces}

Manifolds modeled on the geometry ${\rm Sol}$ are either torus or Klein bottle
bundles over $S^1$ or are composed of two twisted $I$-bundles over
the torus or the Klein bottle. This classification
follows from work of Hempel and Jaco \cite{HempelJaco1972};
see also \cite{Scott1983}.  
Triangulations of torus bundles are discussed in \cite{Lutz2003hpre}.


\section{Hyperbolic Spaces}

According to Thurston \cite{Thurston1982}, ``most'' of the prime $3$-manifolds 
are hyperbolic. In fact, if $L\subset M^3$ is a link in a three-manifold $M^3$
such that $M^3\backslash L$ has a hyperbolic structure, then almost every manifold 
that results from $M^3$ by Dehn surgery along the link $L$ has a hyperbolic structure. 

Hodgson and Weeks \cite{HodgsonWeeks1994} gave a census of 11,031 orientable 
closed hyperbolic $3$-manifolds obtained by Dehn filling on the examples
of a census of cusped hyperbolic three-manifolds by Hildebrand and Weeks 
\cite{HildebrandWeeks1989}. 
The Hildebrand-Weeks census comprises 415 cusped hyperbolic three-manifolds
obtained from gluing five or fewer ideal hyperbolic tetrahedra, a concept
introduced by Thurston (\cite{Thurston1986}, \cite[Ch.~4]{Thurston2002}),
and has later been extended to a census of 6075 cusped hyperbolic 
three-manifolds with up to seven ideal tetrahedra by 
Callahan, Hildebrand, and Weeks \cite{CallahanHildebrandWeeks1999}.

Hyperbolic $3$-manifolds that arise by identifying the faces of a 
Platonic solid were classified by Everitt \cite{Everitt2001pre}.
The best known of these examples is the hyperbolic dodecahedral space
of Weber and Seifert \cite{WeberSeifert1933}.
This space has homology
$({\mathbb Z},{\mathbb Z}_5^{\,3},0,{\mathbb Z})$,
and its smallest known triangulation has $f$-vector
$f=(22,196,348,174)$.


\section{Complex Hypersurfaces and Varieties}
\label{sec:hypersurfaces}

One fruitful and extensively studied way to obtain examples of manifolds 
is via complex hypersurfaces or varieties. In this section we give some
basic facts and examples following Brieskorn (\cite{Brieskorn1966b},
\cite{Brieskorn1966a}), Hirzebruch (\cite{Hirzebruch1966_67},
\cite{HirzebruchMayer1968}), and Milnor (\cite{Milnor1968}, \cite{Milnor1975b}).
(For an introduction into algebraic geometry see Cox, Little, and
O'Shea \cite{CoxLittleOShea1992}.)

For every non-constant complex polynomial $f(z_1,\ldots,z_n)$ in $n$
variables the set 
$$V_f=\{\, z\in {\mathbb C}^n\mid f(z)= 0\,\}$$
is called a \emph{(affine) complex hypersurface}. Let 
$$S^{2n-1}_{(\varepsilon,z^*)}=\{\, z\in {\mathbb C}^n\mid |z_1-z_1^*|^2+\ldots+|z_n-z_n^*|=\varepsilon^2\,\}$$ 
be the sphere of real dimension $2n-1$ in ${\mathbb C}^n$ with center
$z^*$ and small radius~$\varepsilon$. 
(We simply write $S^{2n-1}$ if\, $\varepsilon=1$\, and $z^*$ is the origin.)

A point $z^*\in {\mathbb C}^n$ is \emph{regular}
if not all of the partial derivatives $\partial f/\partial z_j$ vanish at~$z^*$.
In this case, the intersection
$$K=V_f\cap\, S^{2n-1}_{(\varepsilon,z^*)}$$
is a smooth $(2n-3)$-manifold diffeomorphic to the standard
$(2n-3)$-sphere.

If all the partial derivatives $\partial f/\partial z_j$ are zero at $z^*\in{\mathbb C}^n$,
then $z^*$ is a \emph{critical point}, and $z^*$ is called \emph{isolated}
if there are no other critical points in a small neighborhood of $z^*$.

A subset $V\subset {\mathbb C}^n$ is a \emph{(affine) variety} or a \emph{complex algebraic set}
if it arises as an intersection $V=V_{f_1}\cap\ldots\cap V_{f_m}$
of complex hypersurfaces corresponding to a collection $\{ f_1,\ldots f_m\}$
of polynomial functions on ${\mathbb C}^n$.
A point $z^*\in {\mathbb C}^n$ is \emph{non-singular} 
or \emph{simple} if the matrix $(\partial f_i/\partial x_j)$
has maximal rank at $z^*$, and \emph{singular} otherwise.

If $V$ and $V'$ are two algebraic sets in ${\mathbb C}^n$, 
then so is their union $V\cup V'$. \linebreak
A non-empty algebraic set $V$ is \emph{irreducible} 
if it cannot be expressed as the union of two proper algebraic subsets.

\begin{thm}{\rm (\cite[2.9]{Milnor1968})}
For every isolated singular point $z^*$ of an irreducible variety $V$ 
and sufficiently small $\varepsilon$, the
set\, $K=V\cap\, S^{2n-1}_{(\varepsilon,z^*)}$\, is a smooth $(2n-3)$-mani\-fold 
(possibly the empty set).
\end{thm}
Brauner \cite{Brauner1928} considered the neighborhood boundaries $K$ 
of (isolated) singular points for analyzing polynomials in two complex variables.
For example, if $f(z_1,z_2)=z_1^p+z_2^q$, 
then $K$ is a torus knot of type $(p,q)$ in the $3$-sphere $S^3_{(\varepsilon,0)}$.
Higher-dimensional mani\-folds that arise this way 
were first studied by Brieskorn \cite{Brieskorn1966a}.
He proved that for even $n\geq 2$ the manifold
$$\Sigma^{2n-3}(2,2,\ldots,2,3)\,=\,\{\, (z_1,\ldots,z_n)\in {\mathbb C}^n\mid z_1^2+\ldots +z_{n-1}^2+z_n^3= 0\,\}\,\cap\, S^{2n-1}$$
is homeomorphic to the sphere $S^{2n-3}$.
The result of Brieskorn has the following generalization.
\begin{thm}{\rm (\cite[5.2, 8.1]{Milnor1968})}
If the origin is an isolated critical point of the polynomial 
$f(z_1,\ldots,z_n)$, $n\geq 2$, then the $(2n-3)$-dimensional
manifold $K=f^{-1}(0)\cap S^{2n-1}_{(\varepsilon,0)}$
is $(n-3)$-connected. In particular, $K$ is non-empty for $n\geq 2$,
connected for $n\geq 3$, and simply connected for $n\geq 4$. 
Thus, for $n\neq 3$ the manifold $K$ is
homeomorphic to the sphere $S^{2n-3}$ if and only if $K$ has the
homology of a sphere. 
\end{thm}

However, Hirzebruch \cite[11.3]{HirzebruchMayer1968}
observed that for even $n\geq 4$ and odd $r\equiv\pm3\!\mod 8$ the manifold
$$W^{2n-3}(r)\,=\,\{\, (z_1,\ldots,z_n)\in {\mathbb C}^n\mid z_1^2+\ldots +z_{n-1}^2+z_n^r= 0\,\}\,\cap\, S^{2n-1}$$
is diffeomorphic to Kervaire's sphere, 
which is is an exotic sphere if $r$ and $n$ are not powers of~$2$.
In particular, for $n=6$, the $9$-dimensional Kervaire sphere 
$\Sigma (2,2,2,2,2,3)$ is exotic \cite{Kervaire1960}. For even $n\geq 4$ and 
odd $r\equiv\pm1\!\mod 8$ the manifold $W^{2n-3}(r)$ is diffeomorphic
to the standard $(2n-3)$-sphere. 

The manifold $W^3(3)$ is homeomorphic
to the lens space $L(3,1)$ \cite{Hirzebruch1962_63} and coincides with the 
Brieskorn manifold $M(2,2,3)$; see below. Further special cases are:
$W^{(2n-3)}(0)$ is diffeomorphic to $S^{n-1}\times S^{n-2}$, $W^{2n-3}(1)$
is diffeomorphic to $S^{2n-3}$, and $W^{2n-3}(2)$ is the Stiefel manifold
$V_{n,2}$; see \cite[\S 6]{HirzebruchMayer1968}.

Based on Pham's work \cite{Pham1965} on hypersurfaces of the form
$$\Xi_a(t)=\{\, z\in {\mathbb C}^n\mid z_1^{\,a_1}+\ldots+z_n^{\,a_n}= t\,\},\quad\mbox{with}\quad t\in\mathbb{C},$$
Brieskorn \cite{Brieskorn1966b} studied the manifolds
$$\Sigma (a_1,\ldots,a_n)=\Xi_a(0)\cap \, S^{2n-1},$$
where $a=(a_1,\ldots,a_n)$ is a collection of integers with $a_i\geq 2$.
The space $\Sigma (a_1,\ldots,a_n)$ is a topological sphere if and only if
$\Xi_a(0)$ is a topological manifold, and for $n\geq 4$ Brieskorn provided
a necessary and sufficient condition for this to be the case.
In particular, for $m\geq 2$ and $n=2m+1$
the manifolds $\Sigma (2,\ldots,2,3,6k-1)$ with $k=1,\ldots,|bP_{4m-1}|$ are $(4m-1)$-spheres
and represent all $|bP_{4m-1}|$ different differentiable structures in the cyclic
group $bP_{4m-1}$ of $h$-cobordism classes of $(4m-1)$-spheres; cf.\ \cite{KervaireMilnor1963}.
Thus, for $m=2$ and $m=3$ the spheres $\Sigma (2,2,2,3,6k-1)$ and
$\Sigma (2,2,2,2,2,3,6k-1)$ yield
all $28$ respectively $992$ different classes of exotic $7$- and
exotic $11$-spheres.

\bigskip

Let us return to $3$-dimensional manifolds. According to Milnor
\cite{Milnor1975b}, the geometric type of a $3$-dimensional 
Brieskorn manifold $M(p,q,r)$, $p,q,r\geq 2$, 
which results from the intersection of the complex algebraic Pham-Brieskorn hypersurface
$\Xi_{(p,q,r)}(0)=\,\{\, z\in {\mathbb C}^3\mid z_1^p+z_2^q+z_3^r = 0\,\}$
with the sphere $S^5$ depends only on the sign of the
rational number\, $\frac{1}{p}+\frac{1}{q}+\frac{1}{r}-1$.

For\, $\frac{1}{p}+\frac{1}{q}+\frac{1}{r} > 1$ the $3$-manifold $M(p,q,r)$ is spherical
and its fundamental group has order\, $\frac{4}{pqr}(\frac{1}{p}+\frac{1}{q}+\frac{1}{r})^{-2}$.
These manifolds were first discussed by Brieskorn \cite{Brieskorn1968}
and yield the following examples (see also \cite{Milnor1975b} and \cite{Wagreich1972}).
$M(2,2,r)$ are the lens spaces $L(r,1)$, $M(2,3,3)$ is the prism space $P(2)$,
$M(2,3,4)$ is the octahedral space $S^3/T^*$, and $M(2,3,5)$ is the
spherical dodecahedral space. In the case that we allow $r=0,1$ for $M(2,2,r)$,
then the resulting manifolds are $M(2,2,1)\cong L(1,1)\cong S^3$ and
$M(2,2,0)\cong L(0,1)\cong S^2\times S^1$;
see Hirzebruch \cite[5.7]{HirzebruchMayer1968}.

For\, $\frac{1}{p}+\frac{1}{q}+\frac{1}{r} = 1$ Orlik \cite{Orlik1970} proved that
$M(p,q,r)$ belongs to the geometry ${\rm Nil}$. In particular, the three
possible examples 
$M(2,3,6)$, $M(2,4,4)$, and $M(3,3,3)$ are just the $T^2$-bundles 
$M_1$, $M_2$, and $M_3$ over $S^1$ from Section~\ref{sec:nil_spaces};
see \cite{Milnor1975b}.

Finally, if $\frac{1}{p}+\frac{1}{q}+\frac{1}{r} < 1$, then $M(p,q,r)$ is
modeled on $\widetilde{SL}(2,{\mathbb R})$. For a detailed account of
this case see \cite{Milnor1975b}. Not all Brieskorn manifolds
of type $\widetilde{SL}(2,{\mathbb R})$ are distinct:
For example, $M(2,9,18)\cong M(3,5,15)$.

\bigskip

The spherical Brieskorn manifolds yield all spherical $3$-manifolds
that arise as quotients $S^3/G$ with respect to finite subgroups $G$ of $S^3$, 
except for the truncated cube space $S^3/O^*$ and the prism spaces $P(k)$, $k\geq 3$.
However, these manifolds have simple descriptions
as complex hypersurfaces (cf., \cite[p.~80]{Milnor1968}). 
The truncated cube space $S^3/O^*$ 
is defined by the equation 
$$z_1^2+z_2^3+z_2z_3^3=0,$$
and the prism spaces $P(k)$ result from the hypersurfaces
$$z_1^2+z_2^2z_3+z_3^{k+1}=0.$$

\bigskip

A generalization of Milnor's result on Brieskorn manifolds
is due to Neumann \cite{Neumann1977}. Let $(\alpha_1,\ldots,\alpha_n)$ be
an $n$-tuple of integers with $\alpha_i\geq 2$
and \linebreak $A=(a_{ij})_{1\leq i\leq n-2,1\leq j\leq n}$
be a matrix of complex numbers such that every $(n-2)\times (n-2)$ subdeterminant is non-zero.
Then the complex algebraic set
$$V_A(\alpha_1,\ldots,\alpha_n)=\{\, z\in {\mathbb C}^n\mid
a_{i1}z_1^{\,\alpha_1}+\ldots+a_{in}z_n^{\,\alpha_n} = 0,\quad i=1,\ldots,n-2\,\}$$
is called a \emph{Brieskorn complete intersection} or a
\emph{generalized Brieskorn variety}, and the manifold
$M^3(\alpha_1,\ldots,\alpha_n)=V_A(\alpha_1,\ldots,\alpha_n)\cap S^{2n-1}$ is 
a \emph{generalized Brieskorn manifold}. In fact,
$M^3(\alpha_1,\ldots,\alpha_n)$ is a $3$-manifold and is independent
of the choice of the matrix $A$; cf.\ Hamm (\cite{Hamm1971}, \cite{Hamm1972}).
For $n>3$, Neumann proved in \cite{Neumann1977} that
$M^3(\alpha_1,\ldots,\alpha_n)$ is an $\widetilde{SL}(2,{\mathbb R})$-manifold
whenever $\sum_{i=1}^n\frac{1}{\alpha_i}<n-2$. 
The case $\sum_{i=1}^n\frac{1}{\alpha_i}>n-2$
cannot occur for $n>3$, and $\sum_{i=1}^n\frac{1}{\alpha_i}=n-2$
is only possible for $M^3(2,2,2,2)$. 
The latter manifold is the ${\rm Nil}$ $T^2$-bundle $M_4$ over $S^1$;
see Section~\ref{sec:nil_spaces}.

\section{Homology $3$-Spheres}
\label{sec:homology_spheres}

The best known class of homology $3$-spheres
are the \emph{Seifert homology spheres}, i.e.,
the homology $3$-spheres that can be Seifert fibered.
Seifert \cite{Seifert1933} showed that such a 
homology sphere necessarily has the $2$-sphere as its orbit surface 
and is therefore of type $\{ Oo,0\mid b;(\alpha_1,\beta_1),\ldots,(\alpha_r,\beta_r)\}$.
These manifolds are homology spheres if and only if the coefficients satisfy the equation
$$b\alpha_1\cdots\alpha_r+\beta_1\alpha_2\cdots\alpha_r
+\alpha_1\beta_2\cdots\alpha_r+\alpha_1\alpha_2\cdots\beta_r=\pm 1.$$
In fact, if the left hand side is equal to $-1$, then reversing the orientation
of the manifolds yields $+1$ for the new left hand side. Therefore,
without loss of generality, we can require the left hand side to be
equal to $+1$. If two of the $\alpha_i$ have a common divisor, say $d$,
then the left hand side is divisible by $d\neq 1$, which is impossible.
Therefore, the $\alpha_i$ have to be pairwise coprime, and then
${\rm gcd}(\alpha_1\cdots\alpha_r,
           \alpha_2\cdots\alpha_r,
           \alpha_1\alpha_3\cdots\alpha_r,\ldots,
           \alpha_1\cdots\alpha_{r-1})=1$,
from which it follows that there are integers $b$, $\beta_1$, \ldots, $\beta_r$
which satisfy the above equation. There is even a solution with $b=0$, since 
\mbox{${\rm gcd}(\alpha_2\cdots\alpha_r,\alpha_1\alpha_3\cdots\alpha_r,\ldots,\alpha_1\cdots\alpha_{r-1})\!=\!1$}.
The $\beta_i$ can be obtained as iterated cofactors by computing the ${\rm gcd}$
with the Euclidean algorithm. Some of the $\beta_i$ will be negative 
and hence not satisfy the normalization $0<\beta_i<\alpha_i$.
This can be adjusted by simultaneously replacing $\beta_i$ with $\beta_i+x_i\alpha_i$
and $b$ by $b-x_i$. It also follows that $\beta_i$ and $\alpha_i$ are
coprime for every $i=1,\ldots,r$.

\begin{thm}{\rm (Seifert~\cite{Seifert1933})}
For any set of $r\geq 3$ pairwise coprime integers\, $\alpha_1,\ldots,\alpha_r\geq 2$,
there is (up to reversing the orientation) a unique
\textbf{Seifert fibered homology sphere}\, $\Sigma(\alpha_1,\ldots,\alpha_r)$\, 
with $r$ exceptional fibers of the given multiplicities, 
and every fibered homology sphere arises this way. 
The only fibered homology sphere 
with finite fundamental group is the Poincar\'e sphere\, $\Sigma(2,3,5)$.
\end{thm}

All fibered homology spheres with infinite fundamental 
group are modeled on the geometry $\widetilde{SL}(2,{\mathbb R})$.
If $r=3$ and $\alpha_1$, $\alpha_2$, and $\alpha_3$ are coprime integers, then
the fibered homology sphere\, $\Sigma(\alpha_1,\alpha_2,\alpha_3)$\, is homeomorphic 
to the Brieskorn manifold $M(\alpha_1,\alpha_2,\alpha_3$)
and is called a \emph{Brieskorn homology sphere}.
Neumann and Raymond \cite{NeumannRaymond1978} 
showed that every Seifert homology sphere $\Sigma(\alpha_1,\ldots,\alpha_r)$
is diffeomorphic to the generalized Brieskorn manifold with the same
parameters.

\bigskip

Besides Seifert fibered homology spheres, there is an abundance of
hyperbolic homology spheres. As mentioned above, standard techniques for 
constructing arbitrary homology spheres
are Dehn surgery on a link in the $3$-sphere $S^3$ 
(see Dehn \cite{Dehn1910}, Seifert and Threlfall \cite[\S 65]{SeifertThrelfall1934},
Lickorish \cite{Lickorish1962} and Wallace \cite{Wallace1960})
or three-fold branched coverings of $S^3$ over a knot 
(see Alexander \cite{Alexander1919-1920}, 
Hilden \cite{Hilden1974}, Montesinos \cite{Montesinos1976},
and Izmestiev and Joswig \cite{IzmestievJoswig2001pre}).
The connection between knots, links, and $3$-manifolds
attracted a lot of interest. See the books
of Rolfsen \cite{Rolfsen1976}, Prasolov and Sossinsky \cite{PrasolovSossinsky1997},
and Saveliev \cite{Saveliev1999}. The introduction of various new
invariants for $3$-manifolds in the last twenty years 
(see, in particular, the survey of Fintushel and Stern \cite{FintushelStern1990a}) 
has led to rapid progress in low-dimensional topology.

\bibliography{}

\bigskip
\medskip

\noindent
Frank H.\ Lutz\\
Technische Universit\"at Berlin\\
Fakult\"at II - Mathematik und Naturwissenschaften\\
Institut f\"ur Mathematik, Sekr. MA 6-2\\
Stra\ss e des 17.\ Juni 136\\
D-10623 Berlin\\
{\tt lutz@math.tu-berlin.de}

\end{document}